\documentclass[12pt,a4paper]{article}
\usepackage{amsmath}
\usepackage{amssymb}
\usepackage{amsthm}
\usepackage{amsfonts}
\usepackage{latexsym}

\theoremstyle{plain}
\newtheorem{thm}{Theorem}[section]

\newtheorem{lem}{Lemma}[section]

\theoremstyle{definition}
\newtheorem{defn}{Definition}[section]

\theoremstyle{remark}

\title{Local trace formulae and scaling asymptotics\\ in Toeplitz quantization, II}
\author{Roberto Paoletti\footnote{\noindent{\bf Address:}
Dipartimento di Matematica e Applicazioni, Universit\`a degli Studi
di Milano Bicocca, Via R. Cozzi 53, 20125 Milano,
Italy; {\bf e-mail}: roberto.paoletti@unimib.it }}
\date{}

\begin{document}

\maketitle

\begin{abstract}
In the spectral theory of positive elliptic
operators, an important role is played by certain smoothing kernels
related to the Fourier transform of the trace of
a wave operator, which may be heuristically interpreted as
smoothed spectral projectors asymptotically drifting to the right of
the spectrum. In the setting of Toeplitz quantization, we consider
analogues of these, where the wave operator is replaced by the Hardy space
compression of a linearized Hamiltonian flow, possibly composed with a family of
zeroth order Toeplitz operators. We study the local asymptotics of these
smoothing kernels, and specifically how they concentrate on the fixed
loci of the linearized dynamics.
\end{abstract}

\section{Introduction}

The object of this paper are the local asymptotics of certain smoothing kernels in
Berezin-Toeplitz geometric quantization \cite{ber}, \cite{bg}, \cite{z-id}; these may be seen
as local contributions to the singularities of the distributional traces of certain
\lq wave type\rq\, operators, and these singularities are related to the dynamics
of a linearized flow.

In the spectral theory of a positive elliptic first-order self-adjoint pseudodifferential operator $P$ on a compact
manifold, one is led to consider
smoothing operators of the form
\begin{equation}
\label{eqn:smoothed-basic-kernel}
S_{\chi\,e^{-i\lambda(\cdot)}}=:\int_{-\infty}^{+\infty}\chi(\tau)\,e^{-i\lambda\tau}U(\tau)\,d\tau,
\end{equation}
where in this case
$U(\tau)=:e^{i\tau P}$, and $\chi$ is a compactly supported smooth function on $\mathbb{R}$
\cite{h}, \cite{dg}, \cite{gs}. It is suggestive to think of
(\ref{eqn:smoothed-basic-kernel}) as a \lq moving smoothed spectral projector\rq,\, drifting to the right
of the spectrum as $\lambda\rightarrow +\infty$.

If $\lambda_1\le \lambda_2\le\cdots$ are the eigenvalues of $P$, the distributional trace
$\mathrm{tr}\big(U(\tau)\big)=:\sum _je^{i\lambda_j\tau}$ is a well-defined
tempered distribution on $\mathbb{R}$; as a smooth function of $\lambda$,
the trace of (\ref{eqn:smoothed-basic-kernel}) is
the Fourier transform of the compactly supported distribution $\chi\cdot \mathrm{tr}(U)$.
Therefore,
its asymptotics as $\lambda\rightarrow \infty$ yield information about
the singularities of $\mathrm{tr}(U)$ on the support of $\chi$,
hence about the asymptotic distribution of the eigenvalues \cite{c}, \cite{dg}, \cite{col}.

Here we shall study certain analogues of these
operators in Toeplitz quantization.
In this setting, one typically considers a connected complex d-dimensional projective manifold $M$ with
an Hermitian ample line bundle $(B,h)$, such that the unique compatible covariant derivative $\nabla$
on $B$ has curvature $\Theta=-2i\,\omega$, where $\omega$ is a K\"{a}hler form
(one can generalize this picture to almost complex compact symplectic manifolds \cite{bg}, \cite{sz},
but for ease of exposition we shall confine our discussion to the complex projective category).
The symplectic manifold
$(M,\omega)$ is viewed as a classical phase space, and the spaces $H^0\left(M,B^{\otimes k}\right)$
of global holomorphic sections of tensor powers of $B$, with the naturally induced Hilbert structure,
as associated quantum spaces at Planck's constant $\hbar=1/k$.

A \lq classical observable\rq\, is a smooth
real function $f$ on $M$, and a \lq quantum observable\rq\, is a family of self-adjoint operators
$T^{(k)}:H^0\left(M,B^{\otimes k}\right)\rightarrow H^0\left(M,B ^{\otimes k}\right)$.
The Berezin-Toeplitz quantization $T_f^{(k)}:H^0\left(M,B^{\otimes k}\right)\rightarrow H^0\left(M,B^{\otimes k}\right)$
of $f$ is the zeroth order Toeplitz operator
given by the compression of the multiplication operator
by $f$ in $L^2\left(M,B^{\otimes k}\right)$ (the Hilbert space of square summable sections of $B^{\otimes k}$)
with the orthogonal projector onto $H^0\left(M,B^{\otimes k}\right)$.
One is then led to investigate how
the Hamiltonian dynamics of $f$ captures in various ways
the asymptotics of the corresponding quantum dynamics.

%For example, one may consider the asymptotics of
%$T^{(k)}_f:H^0\left(M,A^{\otimes k}\right)\rightarrow H^0\left(M,A^{\otimes k}\right)$ in the semiclassical limit,
%that is, as $k\rightarrow +\infty$; or else, one may take the Hilbert space direct sum
%$H(M,A)=:\bigoplus _{k\ge 0}H^0\left(M,A^{\otimes k}\right)$ and associate to $f$ a \lq collective quantum observable\rq\,
%on $H(M,A)$, and study the asymptotic spectral properties of the latter.

%Let $\widetilde{T}_f:H(M,A)\rightarrow H(M,A)$
%be the first-order Toeplitz operator defined by $\widetilde{T}_f^{(k)}=:k\,T^{(k)}_f$ on $H^0\left(M,A^{\otimes k}\right)$.
%If $f>0$ and $\lambda_1\le \lambda_2\le \cdots$ are the eigevalues of $\widetilde{T}_f$, then
%$\lambda_j\uparrow +\infty$ as $j\rightarrow +\infty$; hence $\sum_je^{i\lambda_j\tau}$ is a well-defined
%tempered distribution on $\mathbb{R}$. Its singularities relate to the asymptotic distribution of the
%spectrum for $j\rightarrow +\infty$, and are
%encapsulated in
%the trace formula of \cite{bg} in terms of Poincar\'{e} map data (in the general setting of so-called Toeplitz structures,
%and in the framework of Fourier-Hermite distributions). The analogue of (\ref{eqn:smoothed-basic-kernel})
%in this setting is obtained with $U(\tau)=e^{i\tau \widetilde{T}_f}$; here one first quantizes $f$
%to $\widetilde{T}_f$, and then exponentiates the latter to a family of unitary Fourier-Toeplitz operators.

For example, one wishes to quantize the Hamiltonian flow $\phi_\tau^M:M\rightarrow M$ of $f$.
While there is a natural lift $\phi^B_\tau:B\rightarrow B$ preserving $h$, the induced unitary action $V(\tau)$
on $L^2\left(M,B^{\otimes k}\right)$
doesn't generally preserve the closed subspace
$H^0\left(M,B^{\otimes k}\right)$. A natural procedure
to circumvent this obstruction was proposed in \cite{z-id}, as follows.

Consider the Hilbert space direct sums $L(M,B)=:\bigoplus_{k\ge 0}L^2\left(M,B^{\otimes k}\right)$ and
$H(M,B)=:\bigoplus_{k\ge 0}H^0\left(M,B^{\otimes k}\right)$, and let $P:L(M,B)\rightarrow H(M,B)$ be the orthogonal projection.
By the above, the compression $P\circ V(\tau)\circ P$ is generally non-unitary as an endomorphism of $H(M,B)$.
Nonetheless, there is a canonically constructed smooth family of $S^1$-invariant
zeroth order Toeplitz operators $R(\tau)$
such that the composition $R(\tau)\circ P\circ V(\tau)\circ P:H(M,B)\rightarrow H(M,B)$ is (essentially)
unitary.
This motivates considering analogues of
(\ref{eqn:smoothed-basic-kernel}) where $U(\tau)$ is replaced by $R(\tau)\circ P\circ V(\tau)\circ P$,
for a general
smooth family of zeroth order Toeplitz operators $R(\tau)$.

In this paper, we shall describe how the
diagonal asymptotics of the latter smoothing kernels concentrate
on appropriate fixed loci of the linearized dynamics of $f$; we shall not require
$R(\tau)$ to be invariant for the $S^1$-action given by fibrewise scalar
multiplication (so these operators belong to a more general class than
those discussed above).

These pointwise asymptotics
may be viewed
as local contributions to the Fourier transform of the  distributional trace of $U$, which motivates
the heuristic term \lq local trace formula\rq. Philosophically, a local trace formula yields information
about the asymptotic distribution of the eigenfunctions, and not just of the corresponding eigenvalues.
Although we won't make this explicit here,
a genuine trace formula may then be obtained by a global integration, as was done in
Corollary 1.1 of \cite{p-ltf} for the special class of Hamiltonians generating holomorphic flows.

For a more precise discussion, let us now backtrack and lift the analysis to the circle bundle of $B$; this approach follows the
philosophy and the techniques
of \cite{bg}, \cite{z-szego}, \cite{bsz}, \cite{sz}, and has the advantage of replacing line bundle
sections with functions, and making available the microlocal description of the Szeg\"{o} kernel as an FIO
in \cite{bs}.

Thus as above let $M$ a d-dimensional connected complex projective manifold, and let
$(B,h)$ be a positive Hermitian line bundle on it.
If $B^\vee$ is the dual line bundle, let $X\subseteq B^\vee$ be the unit circle bundle for the induced
Hermitian structure, with projection $\pi:X\rightarrow M$. Then $X$ is a strictly pseudoconvex domain,
and the connection form $\alpha$ is a
contact form on it; in particular, $dV_M=:(1/\mathrm{d}!)\,\omega^{\wedge \mathrm{d}}$
is a volume form on $M$, and $d\mu_X=:(1/2\pi)\,\alpha\wedge \pi^*(dV_M)$ is a volume form on $X$.

Since the Hardy space $H(X)\subseteq L^2(X)$ is naturally unitarily isomorphic to $H(M,B)$,
Berezin-Toeplitz operators on $H(M,B)$ may be identified with a class of Toeplitz operators on $X$
in the sense of \cite{bg}. These are operators of the form $\Pi\circ Q\circ \Pi$, where
$Q$ is a pseudodifferential operator of classical type on $X$,
and $\Pi:L^2(X)\rightarrow H(X)$ is the orthogonal
projector; $\Pi$ is called the \textit{Szeg\"{o} projector} of $X$, and its distributional kernel
$\Pi\in \mathcal{D}'(X\times X)$ - with abuse of notation -
the \textit{Szeg\"{o} kernel} of $X$.
%Then $T_f$ and $\widetilde{T}_f$ correspond respectively to
%$\mathcal{T}_f=:\Pi\circ M_f\circ \Pi$, and $\widetilde{\mathcal{T}}_f=:D_\theta\circ \mathcal{T}_f$;
%here $\partial/\partial\theta$ is the generator of the structure $S^1$-action, and $D_\theta=:-i\partial/\partial\theta$
%is the \textit{number operator} (it is $k\cdot\mathrm{id}$ on the $k$-th isotype $H_k(X)$).
%In particular, $\mathcal{T}_f$ and $\widetilde{\mathcal{T}}_f$ are a zeroth order and a first order
%self-adjoint Toeplitz operator on $X$, respectively.

Any real $f\in \mathcal{C}^\infty(M)$ generates a 1-parameter group $\phi^M_\tau$ ($\tau\in \mathbb{R}$)
of Hamiltonian symplectomorphisms of $M$, which lifts canonically to a 1-parameter group $\phi^X_\tau$
of contactomorphisms of $X$. Infinitesimally, if $\upsilon_f$ is the Hamiltonian vector field of $f$
with respect to $2\omega$, and if $\upsilon_f^\sharp$ is its horizontal lift to $X$, then
\begin{equation}
\label{eqn:contact-vector-field}
\widetilde{\upsilon}_f=:\upsilon_f^\sharp-f\cdot \partial/\partial\theta
\end{equation}
is a contact vector field on
$(X,\alpha)$ lifting $\upsilon_f$; thus $\widetilde{\upsilon}_f$ and its flow $\phi^X$
depend on $f$, and not merely on $\upsilon_f$ and $\phi^M$ (changing $f$ by a constant will modify the linearization).

The operators $P\circ V(\tau)\circ P:H(M,B)\rightarrow H(M,B)$ correspond under the
isomorphism $H(X)\cong H(M,B)$ to the compressed pull-backs of the contact flow,
$T_\tau=:\Pi\circ \left(\phi^X_{-\tau}\right)^*\circ \Pi:H(M,B)\rightarrow H(M,B)$.
Motivated by the theory in \cite{z-id},
we consider smooth families of operators of the form
\begin{equation}
\label{eqn:smooth-general-family}
U_\tau=:R_\tau\circ T_\tau,
\end{equation}
where $R_\tau$ is a smooth family of zeroth order Toeplitz operators in the sense of \cite{bg}.

One may regard the distributional kernel of $U$ as a distribution on $\mathbb{R}\times X\times X$.
The \textit{trace} of $U$ is then defined functorially as a distribution on $\mathbb{R}$ as
\begin{equation}
\label{eqn:trace}
\mathrm{trace}(U)=:p_*\left(\Delta^*(U)\right),
\end{equation}
where $p:\mathbb{R}\times X\rightarrow \mathbb{R}$ is the projection, and
$\Delta:\mathbb{R}\times X\rightarrow \mathbb{R}\times X\times X$ is the diagonal map
$(\lambda,x)\mapsto (\lambda,x,x)$ \cite{bg}.
Explicitly, for any $\gamma\in \mathcal{C}^\infty_0(\mathbb{R})$ the averaged operator
\begin{equation}
\label{eqn:S-gamma}
S_\gamma=:\int_{-\infty}^{+\infty}\gamma(\tau)\,U_\tau\,d\tau
\end{equation}
is smoothing (see below), and
$$
\langle \mathrm{trace}(U),\gamma\rangle =\mathrm{trace}(S_\gamma)=\int_{-\infty}^{+\infty}S_\gamma(x,x)\,d\mu_X(x).
$$
Wave front considerations show that this functorial description is consistent, and that the singular
support of $\mathrm{trace}(U)$ is contained in the set of periods of $\phi^X$:

\begin{defn}
\label{defn:period}
A real number $\tau\in \mathbb{R}$ is a \textit{period} of $\phi^X$ if there exists $x\in X$ such that
$\phi^X_\tau(x)=x$; let $\mathrm{Per}_X\left(f\right)\subseteq \mathbb{R}$ be the set of periods of $\phi^X$.
Similarly, $\tau\in \mathbb{R}$ is a period of $\phi^M$ if there exists $m\in M$ such that
$\phi^M_\tau(m)=m$, and $\mathrm{Per}_M\left(f\right)\subseteq \mathbb{R}$ will denote the set of periods of $\phi^M$.
\end{defn}

Thus, $\mathrm{Per}_X\left(f\right)\subseteq \mathrm{Per}_M\left(f\right)$, and the inclusion is generally proper.
Furthermore, if $\tau\in \mathrm{Per}_X\left(f\right)$ then the fixed locus $\mathrm{Fix}\left(\phi^X_{\tau}\right)
\subseteq X$ of $\phi^X_{\tau}$
is the inverse image in $X$ of the union $\mathrm{Fix}\left(\phi^M_{\tau}\right)'$
of a collection of connected components of the fixed
locus $\mathrm{Fix}\left(\phi^M_{\tau}\right)\subseteq M$ of $\phi^M_\tau$.

If $U=U(\tau)$ is as just described, we shall view the diagonal asymptotics of the smoothing kernel
(\ref{eqn:smoothed-basic-kernel})
as local contributions to the Fourier transform of $\chi\cdot\mathrm{trace}(U)$.
When $\chi$ is compactly supported near an isolated period $\tau_0$ of $\phi^X$,
these asymptotics will specifically detect the singularity of $\mathrm{trace}(U)$
at $\tau_0$.
Geometrically, a closed
orbit in $M$ will contribute non negligibly to the asymptotics of (\ref{eqn:smoothed-basic-kernel}) only if it
lifts to a loop in $X$.
Explicitly,
\begin{equation}
\label{eqn:fourier-transform}
\langle \chi\cdot\mathrm{trace}(U),e^{-i\lambda(\cdot)}\rangle =\mathrm{trace}(S_{\chi\cdot e^{-i\lambda(\cdot)}})=\int_{-\infty}^{+\infty}S_{\chi\cdot e^{-i\lambda(\cdot)}}(x,x)\,d\mu_X(x),
\end{equation}
where as $\lambda\rightarrow +\infty$ the integrand concentrates on the
fixed locus of $\phi^X_{\tau_0}$, and is rapidly decreasing everywhere
for $\lambda\rightarrow -\infty$.

It is really only for a certain class of periods that we shall exhibit explicit local asymptotics:

\begin{defn}
\label{defn:very-clean-M}
$\tau\in \mathrm{Per}_M\left(f\right)$ is \textit{very clean} if
$\mathrm{Fix}\left(\phi^M_{\tau}\right)$ is a \textit{symplectic} submanifold of $M$, and
\begin{equation}
\label{eqn:clean}
  T_m\mathrm{Fix}\left(\phi^M_{\tau}\right)=\ker \left(d_m\phi^M_{\tau}-\mathrm{id}_{T_mM}\right)
  \end{equation}
  for any $m\in \mathrm{Fix}\left(\phi^M_{\tau}\right)$.
\end{defn}

If $\tau$ is a clean period (i.e., (\ref{eqn:clean}) holds), then it is very clean if and only if
in addition for every
$m\in \mathrm{Fix}\left(\phi^M_{\tau}\right)$ we have
  \begin{equation}
\label{eqn:very-clean}
  \ker \left(d_m\phi^M_{\tau}-\mathrm{id}_{T_mM}\right)\cap \mathrm{im}\left(d_m\phi^M_{\tau}-\mathrm{id}_{T_mM}\right)
  =(0)
  \end{equation}
(\S 4 of \cite{dg}).
Now both (\ref{eqn:clean}) and (\ref{eqn:very-clean}) are local conditions on $\mathrm{Fix}\left(\phi^M_{\tau}\right)$, and
on the other hand our analysis involves only those connected components
of $\mathrm{Fix}\left(\phi^M_{\tau}\right)$ dominated by $\mathrm{Fix}\left(\phi^X_{\tau}\right)$.
We shall thus adopt the following:

\begin{defn}
\label{defn:very-clean-M}
$\tau\in \mathrm{Per}_X\left(f\right)$ is \textit{very clean} if $\mathrm{Fix}\left(\phi^X_{\tau}\right)$
 is a submanifold of $X$, and (\ref{eqn:clean}) and (\ref{eqn:very-clean})
are satisfied
$\forall\,m\in\mathrm{Fix}\left(\phi^M_{\tau}\right)'$.
\end{defn}

Obviously, if $\tau\in $ is $\mathrm{Per}_X\left(f\right)$ is very clean as a period of $\phi^M$, then it
very clean as a period of $\phi^X$.

Actually, in order for the local asymptotics in Theorem \ref{thm:main}
below to hold at a given $x\in X$ with $\phi^X_{\tau}(x)=x$, it is only necessary that (\ref{eqn:clean}) and (\ref{eqn:very-clean})
hold at $m=\pi(x)$. However, for brevity we shall require that $\tau$ is very clean as a period of $\phi_X$; they then
hold uniformly over the fixed locus.

The condition is always satisfied if $\phi^M$ is holomorphic.
By a similar principle any $\tau\in \mathrm{Per}_M\left(f\right)$ -
and thus \textit{a fortiori} any $\tau\in \mathrm{Per}_X\left(f\right)$ -
will be very clean if $\phi^M$ is periodic or, more generally, if $\phi^M$ is an immersed subgroup of a compact Lie
group acting symplectically on $M$. Indeed, averaging over $G$ one can then produce in a standard manner an auxiliary
invariant and compatible almost complex structure, and thus prove that every fixed locus is clean and symplectic
\cite{ms}.

The scaling asymptotics in Theorem \ref{thm:main} will be expressed in terms of Heisenberg local coordinates on $X$
centered at a given $x\in \mathrm{Fix}\left(\phi^X_{\tau}\right)$; Heisenberg local coordinates are defined and
discussed in \cite{sz}. A choice of Heisenberg local coordinates centered at $x$ entails a choice of preferred local
coordinates for $M$ centered at $m=\pi(x)$ \cite{sz}, and therefore a unitary isomorphism $T_mM\cong \mathbb{C}^\mathrm{d}$.
If $\zeta$ is Heisenberg local chart centered at $x\in X$ we shall set $x+(\theta,\mathbf{v})=:\zeta(\theta,\mathbf{v})$ and
$x+\mathbf{v}=:\zeta(0,\mathbf{v})$ where these expressions are defined; if $m=\pi(x)$ we shall also set
$m+\mathbf{v}=:\pi\big(x+(\theta,\mathbf{v})\big)$. Under the induced unitary isomorphism $T_mM\cong \mathbb{C}^\mathrm{d}$,
these expressions have a natural interpretation also for $\mathbf{v}\in T_mM$ of suitably small norm.

Given this, since $d_m\phi^M_{-\tau}:T_mM\rightarrow T_mM$ is a real symplectic automorphism, it corresponds to a symplectic $2\mathrm{d}\times
2\mathrm{d}$ symplectic matrix $A$.
Thus if $\tau$ is very clean as  a period of $\phi^X$,
then the symplectic subspace $T_A=:\ker(A-I)\subseteq \mathbb{R}^{2\mathrm{d}}$
corresponds to $T_m\mathrm{Fix}\left(\phi^M_{\tau}\right)$,
and $N_A=:\mathrm{im}(A-I)$ to the fibre at $m$ of the symplectic normal bundle of $\mathrm{Fix}\left(\phi^M_{\tau}\right)$.
In intrinsic notation, under the same assumption will shall set
$$
T_m=:T_m\mathrm{Fix}\left(\phi^M_{\tau}\right),\,\,\,\,\,\,
N_m=:T_m^{\perp_\omega},
$$
that is, $N_m$ is the \textit{symplectic} normal space to $T_m$, and $T_mM=T_m\oplus N_m$
if $\tau$ is very clean.

A further ingredient is required in the statement of Theorem \ref{thm:main}.
In the near diagonal scaling asymptotics for the equivariant components of the Szeg\"{o} kernel,
the off-diagonal exponential decay is controlled (in Heisenberg local coordinates)
by the function $\psi_2:\mathbb{R}^{2\mathrm{d}}\times \mathbb{R}^{2\mathrm{d}} \rightarrow \mathbb{C}$
defined by
\begin{equation}
\label{eqn:psi-2}
\psi_2(\mathbf{v},\mathbf{w})=:-i\,\omega_0(\mathbf{v},\mathbf{w})-\frac 12\,\|\mathbf{v}-\mathbf{w}\|^2\,\,\,\,\,\,\,\,\,(\mathbf{v},
\mathbf{w}\in \mathbb{R}^{2\mathrm{d}}),
\end{equation}
where $\omega_0$ denotes the standard symplectic structure and $\|\cdot\|$ is the Euclidean norm
(thus, $\omega_0$ is represented by $-J_0$,
where $J_0$ is the matrix corresponding to the standard complex structure).

Actually, $\psi_2$
also controls the diagonal asymptotics for $S_{\chi\,e^{-\lambda\cdot}}$ for normal scaled displacements away from a fixed locus
when $f$ is compatible \cite{p-ltf}, but in the present case of a general Hamiltonian it needs to be modified.
For a general symplectic matrix $A$, let us set
\begin{eqnarray}
\label{eqn:defnQA}
Q_A&=:&I+A^tA,\\
\label{eqn:defnFA}
F_A&=:&J_0\left(A^{-1}-I\right),\\
\label{eqn:defnGA}
G_A&=:&A^t\,(A-I).
\end{eqnarray}

\begin{defn} If $A$ is symplectic and $A=OP$, where $O$ and $P$ are orthogonal and symmetric, respectively,
let us define $\Psi_2^A:\mathbb{R}^{2\mathrm{d}}\rightarrow \mathbb{C}$
by $\Psi_2^A(\mathbf{v})=:\mathbf{v}^t\,\mathfrak{P}_A\,\mathbf{v}$, where
\begin{eqnarray}
\label{eqn:psi-2-A}
\mathfrak{P}_A=:-\left(A^t-I\right)\,O Q_A^{-1}O^t(A-I)-i\,\left(G_A^tQ_A^{-1}F_A-A^tJ_0\right).
\end{eqnarray}
\end{defn}

Suppose for example that $A$ is unitary, that is, orthogonal and symplectic; this will be the case if $f$ is compatible, and
$A$ corresponds to $d_m\phi^M_{-\tau}:T_mM\rightarrow T_mM$
in Heisenberg local coordinates
at a fixed point of $\phi^M_{\tau}$; then $P=I$, $A^tA=I$, and
$$
\mathfrak{P}_A=-\frac 12\,\left(A^t-I\right)\,\big[I-iJ_0\big]\,(A-I)+i\,A^tJ_0.
$$
Thus for any $\mathbf{v}\in \mathbb{R}^{2\mathrm{d}}$ if $A$ is unitary we have
$$
\Psi_2^A(\mathbf{v})=-\frac 12\,\|A\mathbf{v}-\mathbf{v}\|^2-i\omega_0(A\mathbf{v},\mathbf{v})=
\psi_2(A\mathbf{v},\mathbf{v}).
$$

Returning to the general case,
suppose now that $\ker(A-I)\cap \mathrm{im}(A-I)=(0)$, as will be the case under the previous identifications
if $\tau$ is a very clean period. Then
$A-I$ restricts to an automorphism of $N_A=\mathrm{im}(A-I)$, and so
$\Re\left(\Psi_2^A\right)$ is negative definite on $N_A$.

Let $\mathrm{dist}_X$ denote the Riemannian distance function on $X\times X$, extended to a distance between subsets in the
usual manner.

\begin{thm}
\label{thm:main}
Suppose $f\in \mathcal{C}^\infty(M)$, $f>0$. Let $\phi^X$ be the 1-parameter group of contactomorphisms
of $(X,\alpha)$ generated by $f$, and let $\tau_0\in \mathbb{R}$ be a very clean isolated period
 of $\phi^X$. Also, let $R=R_\tau$ be a smooth family of zeroth order Toeplitz operators on
$X$, and set $U_\tau=:R_\tau\circ \left(\phi^X_{-\tau}\right)^*\circ \Pi$.
Finally, given $\gamma\in \mathcal{C}^\infty_0(\mathbb{R})$
define $S_{\gamma}$
by (\ref{eqn:S-gamma}) above.
 Then $S_{\gamma}\in \mathcal{C}^\infty(X\times X)$. Furthermore,
there exists $\epsilon>0$ such that for all
$\chi\in \mathcal{C}^\infty_0\big((\tau_0-\epsilon,\tau_0+\epsilon)\big)$ the
following holds.
\begin{enumerate}
  \item $S_{\chi\,e^{-i\lambda(\cdot)}}(x,x)=O\left(\lambda^{-\infty}\right)$ uniformly on $X$ as $\lambda\rightarrow
   -\infty$.
  \item $S_{\chi\,e^{-i\lambda(\cdot)}}(x,x)=O\left(\lambda^{-\infty}\right)$ uniformly for
  $\mathrm{dist}_X\left(x,\mathrm{Fix}\left(\phi^X_{\tau_0}\right)\right)\ge C\,\lambda^{-7/18}$ as $\lambda\rightarrow
   +\infty$.
  \item Uniformly in $x_0\in \mathrm{Fix}\left(\phi^X_{\tau_0}\right)$ and
  $\mathbf{n}\in N_{\pi(x_0)}\subseteq T_{\pi(x_0)}M$ with $\|\mathbf{n}\|\le C \,\lambda^{1/9}$, the following
  asymptotic expansion holds for $\lambda\rightarrow +\infty$:
  \begin{eqnarray}
  \label{eqn:main}
  \lefteqn{S_{\chi\,e^{-i\lambda(\cdot)}}\left(x_0+\frac{\mathbf{n}}{\sqrt{\lambda}},x_0+\frac{\mathbf{n}}{\sqrt{\lambda}}\right)}\\
  &\sim& \rho_{\tau_0}(x_0)\,
\frac{2\pi\,e^{-i\lambda\tau_0}}{f(m_0)^{\mathrm{d}+1}}\cdot\frac{2^\mathrm{d}}{\sqrt{\det(Q_A)}}
\cdot e^{
f(m_0)^{-1}\Psi_2^A\left(\mathbf{n}\right)}\,\chi(\tau_0)\nonumber\\
  && \cdot \left(\frac \lambda\pi\right)^\mathrm{d} \left[1+\sum_{j\ge 1}\lambda^{-j/2}G_j(x_0,\mathbf{n})\right],\nonumber
  \end{eqnarray}
  where $\rho_\tau$ is the symbol of $R_\tau$,
  and the $G_j$'s are polynomials in $\mathbf{n}$, depending smoothly on $x_0$.
  \item Let us write
  $$
  S_{\chi\,e^{-i\lambda(\cdot)}}\left(x_0+\frac{\mathbf{n}}{\sqrt{\lambda}},x_0+\frac{\mathbf{n}}{\sqrt{\lambda}}\right)=
  \mathfrak{E}_\lambda(x_0,\mathbf{n})+\mathfrak{O}_\lambda(x_0,\mathbf{n}),
  $$
  where $\mathfrak{E}_\lambda$ and $\mathfrak{O}_\lambda$ are even and odd functions of $\mathbf{n}$, respectively.
  Then the asymptotic expansion for $\mathfrak{E}_\lambda(x_0,\mathbf{n})$ (respectively, for
  $\mathfrak{O}_\lambda(x_0,\mathbf{n})$)
  as $\lambda\rightarrow +\infty$ is obtained by collecting all integer powers (respectively,
  all fractional powers) of $\lambda$ in (\ref{eqn:main}).
\end{enumerate}
\end{thm}

Notice that if $R$ is $S^1$-invariant (in particular, if it is the identity)
then so is the diagonal restriction of $S_{\chi\,e^{-\lambda\cdot}}$, which is then a function on $M$.
Also, in view of the above considerations (\ref{eqn:main}) expresses an exponential decay of
the scaling asymptotics along directions symplectically orthogonal to the fixed locus.

This result continues and extends those in \cite{p-weyl}  and \cite{p-ltf}.
Specifically, in \cite{p-weyl} the focus was on the singularity at $\tau=0$, and
its use towards a local Weyl law for Toeplitz operators; in \cite{p-ltf}, local
asymptotics were given for general periods, under restrictive assumptions
on the underlying symplectic dynamics (namely, that is preserves the holomorphic
structure).

Finally, we remark that the previous discussion and results could be phrased
in the more general context of compact quantizable almost complex symplectic manifolds,
adopting the approach of \cite{sz}.

\section{Preliminaries and Notation}

Let the cotangent bundle $T^*X$ be endowed with the canonical symplectic structure
$\omega_{\mathrm{can}}=d\mathbf{p}\wedge d\mathbf{q}$, where $\mathbf{q}$ are local coordinates
on $X$ and $\mathbf{p}$ linear coordinates on the fibers.
Then the cotangent lift $\phi^{T^*X}:T^*X\rightarrow T^*X$
of the contact flow $\phi^X_\tau:X\rightarrow X$ is the Hamiltonian flow
generated by the smooth function $H\big((x,\beta)\big)=:-\beta\left(\widetilde{\upsilon}_f(x)\right)$.

Because $\omega$ is symplectic,
$$
\Sigma=:\big\{(x,r\alpha_x):\,x\in X,\,r>0\big\}
$$
is a closed symplectic cone in $T^*X\setminus(0)$ with respect to the standard symplectic
structure. Since $\phi^X$ is a group of
contactomorphisms, $\phi^{T^*X}$
on $T^*X$ leaves $\Sigma$
invariant, and therefore it induces a flow of Hamiltonian symplectomorphisms
$\phi^\Sigma_\tau:\Sigma\rightarrow \Sigma$. The corresponding Hamiltonian on $\Sigma$
is $\widetilde{f}=:r\cdot f$. Explicitly,
$$
\phi^\Sigma_\tau\big((x,r\alpha_x)\big)=\left(   \phi^X_\tau(x),r\,\alpha_{\phi^X_\tau(x)}\right).
$$

If $T=\Pi\circ Q\circ \Pi$ is a Toeplitz operator, its symbol $\sigma_T:\Sigma\rightarrow \mathbb{C}$ is the
restriction to $\Sigma$ of the symbol of $Q$; it is well-defined by the theory of \cite{bg}.
If $T$ has order $r$, $\sigma_T$ is homogenous of degree $r$; in particular, since $R(\tau)$
is zeroth order, its symbol $\rho_\tau =\sigma_{R(\tau)}$
is a $\mathcal{C}^\infty$ function on $X$.

Let us now discuss some useful matrix identities. Let $A$ be a symplectic $2\mathrm{d}\times 2\mathrm{d}$
matrix, and define $Q_A,F_A,G_A$ by (\ref{eqn:defnQA}), (\ref{eqn:defnFA}), (\ref{eqn:defnGA}).
Define the symmetric matrix:
\begin{eqnarray}
\label{eqn:R-A}
\mathfrak{R}_A&=:&-\frac 12\,\left[\left(A^t-I\right)\,(A-I)+F_A^t\,Q_A^{-1}\,F_A-G_A^t\,Q_A^{-1}\,G_A\right]\nonumber\\
&&-i\,\left[G_A^tQ_A^{-1}F_A-A^tJ_0\right].
\end{eqnarray}
Let $\mathfrak{P}_A$ be as in (\ref{eqn:psi-2-A}). Then

\begin{lem}
\label{lem:RA=PA}
$\mathfrak{R}_A=\mathfrak{P}_A$.
\end{lem}

\begin{proof}
It suffices to consider the real parts.
Since $A^tJ_0A=J_0$, we have
\begin{eqnarray*}
F_A=J_0\,\left(A^{-1}-I\right)=J_0\,A^{-1}\,(I-A)=-A^{t}\,J_0\,(A-I).
\end{eqnarray*}
Thus
\begin{eqnarray}
\label{eqn:partereale}
\lefteqn{\left(A^t-I\right)\,(A-I)+F_A^t\,Q_A^{-1}\,F_A-G_A^t\,Q_A^{-1}\,G_A}\\
&=&\left(A^t-I\right)\,\left[I-J_0AQ_A^{-1}A^tJ_0-AQ_A^{-1}A^t\right]\,(A-I).\nonumber
\end{eqnarray}
Let us write $A=OP$, where $O$ and $P$ are orthogonal and symmetric (and symplectic),
respectively. Being unitary, $O$ (and $O^t=O^{-1}$) commutes with $J_0$; on the other hand
$J_0P^k=P^{-k}J_0$ for any integer $k$, so that
$(I+P^2)^{-1}J_0=\left(I+P^{-2}\right)^{-1}J_0$. Using this, we get
\begin{eqnarray}
\label{eqn:primoaddendo}
-J_0AQ_A^{-1}A^tJ_0&=&-J_0OP\,\left(I+P^2\right)^{-1}\,PO^tJ_0\nonumber\\
&=&O\,\left(I+P^2\right)^{-1}O^t.
\end{eqnarray}
Furthermore,
\begin{eqnarray}
\label{eqn:secondoaddendo}
-A\,Q_A^{-1}A^t&=&-OP\,\left(I+P^2\right)^{-1}PO^t\nonumber\\
&=&-OP^2\,\left(I+P^2\right)^{-1}O^t.
\end{eqnarray}

Let us insert (\ref{eqn:primoaddendo}) and (\ref{eqn:secondoaddendo}) in (\ref{eqn:partereale}).
We get
\begin{eqnarray}
\label{partereale1}
\lefteqn{\left(A^t-I\right)\,(A-I)+F_A^t\,Q_A^{-1}\,F_A-G_A^t\,Q_A^{-1}\,G_A}\\
&=&\left(A^t-I\right)\,\left[I+O\,\left(I+P^2\right)^{-1}O^t-OP^2\,\left(I+P^2\right)^{-1}O^t\right]\,(A-I)\nonumber\\
&=&\left(A^t-I\right)\,O\left[I+\left(I+P^2\right)^{-1}-P^2\,\left(I+P^2\right)^{-1}\right]O^t\,(A-I)\nonumber\\
&=&2\left(A^t-I\right)\,O\,\left(I+P^2\right)^{-1}O^t\,(A-I).\nonumber
\end{eqnarray}
\end{proof}

In the following, we shall write $F$, $G$ and $Q$ for $F_A$, $G_A$ and $Q_A$.

\section{Proof of Theorem \ref{thm:main}}

We may assume $\chi(\tau_0)=1$.

By the theory of \cite{bg}, there exists a family $Q_\tau$
of zeroth order pseudodifferential operators on $X$ such that
$[\Pi,Q_\tau]=0$, and $R_\tau=\Pi\circ Q_\tau\circ \Pi$. Thus
\begin{equation}
\label{eqn:factorize-U-tau}
U_\tau=\Pi\circ Q_\tau\circ \left(\phi^X_{-\tau}\right)^*\circ \Pi.
\end{equation}
Now $Q_\tau\circ \left(\phi^X_{-\tau}\right)^*$ is an FIO associated to the conormal
bundle of the graph of $\phi^X_{-\tau}$; hence its kernel can be written microlocally in the form
\begin{equation}\label{eqn:U-tau-FIO}
Q_\tau\circ \left(\phi^X_{-\tau}\right)^*(x,y)=\int_{\mathbb{R}^{2\mathrm{d}+1}}e^{i[\varphi(\tau,x,\eta)-y\cdot \eta]}
b(\tau,x,y,\eta)\,d\eta,
\end{equation}
where $b(\tau,\cdot,\cdot,\cdot)\in S^0_{\mathrm{cl}}$ for every $\tau$, and $\varphi(\tau,\cdot,\cdot)$
is a generating function for $\phi^{T^*X}_{-\tau}$. Thus in local coordinates
$\varphi$ satisfies the Hamilton-Jacobi equation
$\partial \varphi/\partial \tau=-(\partial\varphi/\partial x)\cdot \widetilde{\upsilon}_f$.
If we view (\ref{eqn:U-tau-FIO}) as a distribution on $\mathbb{R}\times X\times X$, its wave front is therefore
\begin{eqnarray}
\label{eqn:wave-front-U-tau}
\mathrm{W F}\left(Q_\tau\circ \left(\phi^X_{-\tau}\right)^*\right)&=&
 \left\{
 \Big((\tau,t),\big(x,\xi\big),
\phi^{T^*X}_{-\tau} \big(x,-\xi\big)\Big):\right.
\\
  &&(x,\xi)\in T^*X\setminus (0),\,\tau\in \mathbb{R}, \,t=-\xi\left(\widetilde{\upsilon}_f(x)\right)\Big\}.
  \nonumber
\end{eqnarray}

By \cite{bs}, $\Pi$ has wave front
$$
\mathrm{WF}(\Pi)=\big\{(x,r\alpha_x,x,-r\alpha_x):\,x\in X,R>0\big\}.
$$
Therefore, in view of (\ref{eqn:factorize-U-tau}) and (\ref{eqn:wave-front-U-tau})
the wave front of
$U\in \mathcal{D}'(\mathbb{R}\times X\times X)$
is
\begin{eqnarray}\label{eqn:wave-front-U-tau-1}
\mathrm{WF}\big (U\big )&=&
 \left\{
 \Big(\big(\tau,r\,f(x)\big),\big(x,r\alpha_x\big),
\left(\phi^X_{-\tau}(x),-r\,\alpha_{\phi^X_{-\tau}(x)}\right)\Big):\right.
\\
  &&x\in X,\,r>0,\,\tau\in \mathbb{R}\Big\},
  \nonumber
\end{eqnarray}
where of course $f(x)=f(m)$ if $m=\pi(x)$. Since $f>0$ by assumption,
$\big(\tau,r\,f(x)\big)$ is never zero as a cotangent vector to $\mathbb{R}$
at $\tau$.

Let $q:\mathbb{R}\times X\times X\rightarrow X\times X$ be the projection.
For any $\beta\in \mathcal{C}^\infty_0(\mathbb{R})$, $\beta\cdot U$ is a
compactly supported distribution on $\mathbb{R}\times X\times X$, and
$$
S_\beta=q_*(\beta\cdot U)\in \mathcal{D}'(X\times X).
$$
Given  (\ref{eqn:wave-front-U-tau-1}),
$\mathrm{WF}(S_\beta)=\emptyset$ by the functorial properties of wave fronts (\cite{h},
\cite{hor-FIO-I},
Proposition 1.3.4 of \cite{d}); hence $S_\beta\in \mathcal{C}^\infty(X\times X)$.

Fix a small $\delta>0$; then we
can find $\epsilon>0$ so small that the following holds.
First  $\mathrm{Per}_X(f)\cap (\tau_0-\epsilon,\tau_0+\epsilon)=\{\tau_0\}$;
furthermore, if
$\mathrm{dist}_X\left(x,\mathrm{Fix}\left(\phi^X_{\tau_0}\right)\right)\ge \delta/2$,
then $\mathrm{dist}_X\left(x,\mathrm{Fix}\left(\phi^X_{\tau}\right)\right)\ge \delta/4$
for any $\tau\in (\tau_0-\epsilon,\tau_0+\epsilon)$.

Now choose $\chi\in \mathcal{C}^\infty_0\big((\tau_0-\epsilon,\tau_0+\epsilon)\big)$.
As the singular support
of $U(\tau)$ is
%$\mathrm{sing.supp.}\big(U(\tau)\big)=
$\mathrm{graph}\left(\phi^X_{-\tau}\right)\subseteq X\times X$,
we conclude that $\chi(\tau)\cdot U(\tau)(x,x)$ is
$\mathcal{C}^\infty$ in $\tau\in \mathbb{R}$ and $x\in X$
such that $\mathrm{dist}_X\left(x,\mathrm{Fix}\left(\phi^X_{\tau_0}\right)\right)\ge \delta/2$;
hence its Fourier transform in $\tau$ is rapidly decreasing. In other words, in the same range
$$
S_{\chi\,e^{-i\lambda\cdot}}(x,x)=O\left(\lambda^{-\infty}\right)
$$
for $\lambda\rightarrow \infty$.

Let us now consider the situation in a tubular neighborhood of $\mathrm{Fix}\left(\phi^X_{\tau_0}\right)$.
We may then assume $x=x_0+\mathbf{n}$, where $x_0\in \mathrm{Fix}\left(\phi^X_{\tau_0}\right)$ and
$\mathbf{n}\in N_{\pi(x_0)}$ has norm $\le\delta$; the expression $x_0+\mathbf{n}$
is interpreted in a system of Heisenberg local coordinates centered at $x_0$ and smoothly varying with $x_0$,
at least locally. In particular, since Heisenberg local coordinates are isometric at the origin,
$\mathrm{dist}_X(x,x_0)< 2\delta$, say.

We have $U_\tau=R_\tau\circ \Pi_\tau$, where $\Pi_\tau=:\left(\phi^X_{-\tau}\right)^*\circ \Pi:L^2(X)\rightarrow L^2(X)$. In terms of
Schwartz kernels,
$$
\Pi_\tau=\left(  \phi^X_{-\tau}\times \mathrm{id}_X\right)^*(\Pi),
$$
and
$\mathrm{WF}(\Pi_\tau)=
%&=&\big\{(x,r\,\alpha_x,y,-r\,\alpha_y):y\in X, \,y=\phi^X_{-\tau}(x),r>0\big\}\nonumber\\
%&=&
\mathrm{graph}\left(\phi^\Sigma_{-\tau}\right)$.
In standard distributional short-hand,
\begin{eqnarray}
\label{eqn:smoothing-kernel-integral}
S_{\chi\,e^{-i\lambda (\cdot)}}(x,x)
&=&\int_{-\infty}^{+\infty}\chi(\tau)\,
e^{-i\lambda\tau}\big(R_\tau\circ \Pi_\tau\big)(x,x)\,d\tau\\
&=& \int_X \int_{-\infty}^{+\infty}\chi(\tau)\,
e^{-i\lambda\tau}R_\tau\left(x,y\right) \Pi_\tau\left(y,x\right)\,d\mu_X(y)\,d\tau.     \nonumber
%&=&\nonumber
\end{eqnarray}

Define
$$
X'=X'(x_0,\delta)=:\Big\{x\in X:\mathrm{dist}_X(x,x_0)<4\delta\Big\},
$$
$$
X''=X''(x_0,\delta)=:\Big\{x\in X:\mathrm{dist}_X(x,x_0)>3\delta\Big\},
$$
and let
$\left\{\varrho',\varrho''\right\}$ be a partition of unity on $X$ subordinate to the open cover $\{X',X''\}$.
Let
$S_{\chi\,e^{-i\lambda (\cdot)}}\left(x,x\right)'$ and $S_{\chi\,e^{-i\lambda (\cdot)}}(x,x)''$
be given by (\ref{eqn:smoothing-kernel-integral}), with the integrand multiplied by
$\varrho'(y)$ and $\varrho''(y)$, respectively.
Thus $S_{\chi\,e^{-i\lambda (\cdot)}}(x,x)=S_{\chi\,e^{-i\lambda (\cdot)}}(x,x)'+
S_{\chi\,e^{-i\lambda (\cdot)}}(x,x)''$.

\begin{lem}
\label{lem:prima-riduzione}
$S_{\chi\,e^{-i\lambda (\cdot)}}(x,x)''=O\left(\lambda^{-\infty}\right)$ as $\lambda\rightarrow \infty$,
uniformly for $x_0\in \mathrm{Fix}\left(\phi^X_{\tau_0}\right)$, $\|\mathbf{n}\|\le\delta$, and $x=x_0+\mathbf{n}$.
\end{lem}

\begin{proof}
As we have remarked, we may assume
$\mathrm{dist}_X(x,x_0)< 2\delta$.
On the other hand,
$\mathrm{dist}_X(y,x_0)>3\delta$ where $\varrho''(y)\neq 0$, whence
$\mathrm{dist}_X\big(y,x\big)\ge \delta$ if $y\in X''$.

Because the singular support of $R_\tau$ is the diagonal
in $X\times X$,
$$g_{x,\tau}(y)=:\varrho''(y)\,R_\tau(x,y)$$ is $\mathcal{C}^\infty$,
uniformly so in the given range.
Define $\Pi^\tau=:\left(\mathrm{id}_X\times \phi^X_{-\tau}\right)^*(\Pi)$
(thus $\Pi^\tau(y,y')=\overline{\Pi_\tau(y',y)}$); as $\Pi$ is regular, so is
$\Pi^\tau$.
We have:
\begin{eqnarray}
\label{eqn:regulare-decade}
S_{\chi\,e^{-i\lambda (\cdot)}}(x,x)''
&=&\int_{-\infty}^{+\infty}e^{-i\lambda\tau}\chi(\tau)\,\left[\int _X \Pi_\tau(y,x)\,g_{x,\tau}(y)\,d\mu_X(y)\right]\,d\tau
\nonumber\\
&=&\int_{-\infty}^{+\infty}e^{-i\lambda\tau}\chi(\tau)\,\overline{\Pi^\tau\left(\overline{g_{x,\tau}}\right)(x)}\,d\tau,
\end{eqnarray}
the Fourier transform in $\tau$ of a smooth function compactly supported in $\tau$.
The statement follows.
\end{proof}

Thus  $S_{\chi\,e^{-i\lambda (\cdot)}}(x,x)\sim S_{\chi\,e^{-i\lambda (\cdot)}}(x,x)'$. To estimate the latter
asymptotically, we may now work in Heisenberg local coordinates centered at $x_0$,
and set $y=x_0+(\theta,\mathbf{v})$, where $(\theta,\mathbf{v})\in (-\pi,\pi)\times \mathbb{C}^\mathrm{d}$,
and $\big\|(\theta,\mathbf{v})\big\|
=O(\delta)$. Suppose
$d\mu_X(y)=\mathcal{V}(\theta,\mathbf{v})\,d\theta\,d\mathbf{v}$ in local coordinates.

By \cite{bs}, $\Pi$ is an FIO of the form
\begin{equation}
\label{eqn:pi-FIO}
\Pi\left(x',x''\right)=\int _0^{+\infty}e^{it\,\psi(x',x'')}\,s\left(t,x',x''\right)\,dt,
\end{equation}
where $\psi$ is a complex phase of positive type, and the amplitude
admits an asymptotic expansion
$s\left(t,x',x''\right)\sim \sum_{j\ge 0}t^{\mathrm{d}-j}\,s_j\left(x',x''\right)$.
Hence microlocally
 \begin{equation}
\label{eqn:R-tau-FIO}
R_\tau\left(x',x''\right)=\int _0^{+\infty}e^{it\,\psi(x',x'')}\,a_\tau\left(t,x',x''\right)\,dt,
\end{equation}
where $a_\tau\left(t,x',x''\right)\sim \sum_{j\ge 0}t^{\mathrm{d}-j}\,a_j\left(x',x''\right)$, with
$a_0\left(x',x'\right)=\rho_\tau\left(x'\right)\,s_0\left(x',x'\right)$. Also,
\begin{equation}
\label{eqn:Pi-tau-FIO}
\Pi_\tau\left(x',x''\right)=\int _0^{+\infty}e^{iu\,\psi_\tau(x',x'')}\,s_\tau\left(u,x',x''\right)\,du,
\end{equation}
where $\psi_\tau\left(x',x''\right)=:\psi\left(\phi^X_{-\tau}\left(x'\right),x''\right)$,
$s_\tau\left(u,x',x''\right)=:s\left(u,\phi^X_{-\tau}\left(x'\right),x''\right)$.

Let us write $\sim$ for \lq has the same asymptotics as\rq.
Working in the neighborhood of $x_0$, we can thus write
\begin{eqnarray}
\label{eqn:S-microlocal-general}
\lefteqn{S_{\chi\,e^{-i\lambda (\cdot)}}(x,x)\sim S_{\chi\,e^{-i\lambda (\cdot)}}(x,x)'}\\
&\sim&\int_0^{+\infty}\int_0^{+\infty}\left(\int_{-\infty}^{+\infty}\int_{-\pi}^\pi\int_{\mathbb{C}^\mathrm{d}}
e^{i\Phi}\,\mathcal{A}\,d\tau\,d\theta\,d\mathbf{v}\right)\,dt\,du,\nonumber
\end{eqnarray}
where
\begin{equation}
\label{eqn:defn-Phi-A}
\begin{array}{ccl}
  \Phi & =: & t\,\psi\big(x,x_0+(\theta,\mathbf{v})\big)
+u\,\psi\left(\phi^X_{-\tau}\big(x_0+(\theta,\mathbf{v})\big),x\right)-\lambda\,\tau,\,\,\,\,
  \mathrm{and}  \\
  \mathcal{A} & =: & a_\tau\big(t,x,x_0+(\theta,\mathbf{v})\big)\,
s_\tau\big(u,x_0+(\theta,\mathbf{v}),x\big)\,\chi(\tau)\,\varrho'(\theta,\mathbf{v})\,
\mathcal{V}(\theta,\mathbf{v}),
\end{array}
\end{equation}
where we have written $\varrho'(\theta,\mathbf{v})$ for $\varrho'\big(x_0+(\theta,\mathbf{v})\big)$.
In particular, integration in $(\tau,\theta,\mathbf{v})$ is compactly supported.

To proceed, we need to dwell on the argument $\phi^X_{-\tau}\big(x_0+(\theta,\mathbf{v})\big)$
appearing in (\ref{eqn:defn-Phi-A}). It is convenient to make the change of variables
$\tau\rightsquigarrow\tau+\tau_0$, where now $|\tau|<\epsilon$.

In Heisenberg local coordinates, $d_{m_0}\phi^M _{-\tau_0}(\mathbf{v})=A\mathbf{v}$.
If $m_0=\pi(x_0)$, for $\mathbf{v}\sim \mathbf{0}$ in $\mathbb{C}^\mathrm{d}$
we have
\begin{eqnarray}
\label{eqn:action-in-heis-local-coordinates-M}
\phi^M _{-\tau_0}(m_0+\mathbf{v})=m_0+\Big(A\mathbf{v}+O\left(\|\mathbf{v}\|^2\right)\Big).
\end{eqnarray}
Lifting to $X$,
\begin{eqnarray}
\label{eqn:action-in-heis-local-coordinates-X}
\phi^X _{-\tau_0}\Big(x_0+(\theta,\mathbf{v})\Big)=x_0+\Big(\theta+\vartheta(\mathbf{v}),
A\mathbf{v}+O\left(\|\mathbf{v}\|^2\right)\Big),
\end{eqnarray}
for a certain $\mathcal{C}^\infty$ real function $\vartheta$ defined near the origin in
$\mathbb{C}^\mathrm{d}$.

\begin{lem}
\label{lem:vartheta-lift}
$\vartheta(\mathbf{v})=O\left(\|\mathbf{v}\|^3\right)$ for $\mathbf{v}\sim \mathbf{0}$.
\end{lem}

\begin{proof}
For $\mathbf{v}\in \mathbb{C}^\mathrm{d}$ of unit norm, consider the smooth path
$\gamma (s)=:m_0+s\,\mathbf{v}$ in $M$, defined for sufficiently small
$s$. Let $\gamma^\sharp$ be the unique horizontal lift of $\gamma$ to $X$
such that $\gamma^\sharp(0)=x_0$. By Lemma 2.4 of \cite{dp},
$\gamma^\sharp(s)=x_0+
\big(g(s\mathbf{v}),s\mathbf{v}\big)$,
where $g(s\mathbf{v})=O\left(s^3\right)$.

Let $\gamma_{-\tau_0}=:\phi^M_{-\tau_0}\circ \gamma$. By (\ref{eqn:action-in-heis-local-coordinates-M}),
$\gamma_{-\tau_0}(s)=m_0+\left(sA\mathbf{v}+O\left(s^2\right)\right)$.
Let $\gamma^\sharp_{-\tau_0}$ be the unique horizontal lift of $\gamma_{-\tau_0}$ to $X$
such that $\gamma^\sharp_{-\tau_0}(0)=x_0$. Thus
\begin{equation}
\label{eqn:action-lift-path}
\gamma^\sharp_{-\tau_0}(s)=x_0+
\Big(g\left(sA\mathbf{v}+O\left(s^2\right)\right),sA\mathbf{v}+O\left(s^2\right)\Big).
\end{equation}

On the other hand, since $\phi^X_{-\tau_0}$ is a contactomorphism covering $\phi^M_{-\tau_0}$ and fixing $x_0$,
the composition $\phi^X_{-\tau_0}\circ \gamma^\sharp$
is also a horizontal lift of $\gamma_{-\tau_0}$, and satisfies
$\phi^X_{-\tau_0}\circ \gamma^\sharp(0)=\phi^X_{-\tau_0}(x_0)=x_0$;
therefore $\phi^X_{-\tau_0}\circ \gamma^\sharp=\gamma^\sharp_{-\tau_0}$.
Given this, applying
(\ref{eqn:action-in-heis-local-coordinates-X}) with $\theta=g(s\mathbf{v})$ we obtain
\begin{eqnarray}
\label{eqn:action-lift-path-1}
\gamma^\sharp_{-\tau_0}(s)&=&\phi^X_{-\tau_0}\Big(x_0+
\big(g(s\mathbf{v}),s\mathbf{v}\big)\Big)\\
&=&\Big(\vartheta(s\mathbf{v})+g(s\mathbf{v}),s\,A\mathbf{v}+O\left(s^2\right)
\Big).\nonumber
\end{eqnarray}
Comparing (\ref{eqn:action-lift-path}) and (\ref{eqn:action-lift-path-1}),
we get $\vartheta(s\mathbf{v})=O\left(s^3\right)$.

\end{proof}

Since Heisenberg local coordinates are horizontal at the origin, by (\ref{eqn:contact-vector-field})
we have
$$
\widetilde{\upsilon}_X(x_0)=\big(-f(m_0),\upsilon_f(m_0)\big)\in \mathbb{R}\times \mathbb{C}^\mathrm{d},
$$
where $\upsilon_f(m_0)$ is identified with its local coordinate expression at $m_0$.
For $(\tau,\theta,\mathbf{v})\sim (0,0,\mathbf{0})\in \mathbb{R}\times \mathbb{C}^\mathrm{d}$,
by Corollary 2.2 of \cite {pao-torus} we get
\begin{eqnarray}
\label{eqn:action-in-tau}
\lefteqn{\phi^X _{-\tau-\tau_0}\Big(x_0+(\theta,\mathbf{v})\Big)=
\phi^X _{-\tau}\Big(x_0+\Big(\theta+\vartheta(\mathbf{v}),
A\mathbf{v}+O\left(\|\mathbf{v}\|^2\right)\Big)\Big)}\nonumber\\
&=&x_0+\Big(\theta+\tau\,f(m_0)+\tau\,\omega_m\big(\upsilon_f(m),A\mathbf{v}\big)+O\left(\|(\tau,\mathbf{v})\|^3\right),
\nonumber\\
&&A\mathbf{v}-\tau\,\upsilon_f(m)+O\left(\|(\tau,\mathbf{v})\|^2\right)\Big).
\end{eqnarray}

Let us consider the asymptotics for $\lambda\rightarrow -\infty$. Given (\ref{eqn:defn-Phi-A}), we have
\begin{equation}
\label{eqn:partial-infty}
\begin{array}{ccl}
  \partial_\tau\Phi & = & u\,\partial_\tau\psi_\tau\big(x_0+(\theta,\mathbf{v}),x\big)+|\lambda|,\\
  \partial_\theta\Phi & = & t\,\partial_\theta\psi\big(x,x_0+(\theta,\mathbf{v})\big)+u\,
  \partial_\theta \psi_\tau\big(x_0+(\theta,\mathbf{v}),x\big).
\end{array}
\end{equation}

Recalling that $x=x_0+\mathbf{n}$ and that
$d_{(x,x)}\psi=(\alpha_x,-\alpha_x)$, in view of (\ref{eqn:action-in-tau})
we have
\begin{equation}
\begin{array}{ccl}
  \partial_\tau\psi_\tau\big(x_0+(\theta,\mathbf{v}),x\big) &= & f(m_0)\cdot\big[1+O(\epsilon,\delta)\big] \\
  \partial_\theta\psi\big(x,x_0+(\theta,\mathbf{v})\big) & =& -1+O(\epsilon,\delta) \\
  \partial_\theta\psi_\tau\big(x_0+(\theta,\mathbf{v}),x\big) & = & 1+O(\epsilon,\delta) ,
\end{array}
\end{equation}
since $\big\|(\theta,\mathbf{n},\mathbf{v})\big\|=O(\delta)$, $|\tau|<\epsilon$.

Therefore,
\begin{eqnarray*}
\nabla_{\theta,\tau}\Phi&=&\Big(t\,\big(-1+O(\tau,\delta)\big)+u\,\big(1+O(\tau,\delta)\big),
u\,f(m_0)\cdot\big[1+O(\tau,\delta)\big]+|\lambda|\Big)^t\\
&=&\big(-t+u,
u\,f(m_0)+|\lambda|\big)^t+\mathbf{R}
\end{eqnarray*}
where $\|\mathbf{R}\|\le \sqrt{t^2+u^2}\cdot O(\tau,\delta)$.
Since $f(m_0),u,t>0$, we have if $\delta$ is sufficiently small
\begin{eqnarray*}
\left\|\big(-t+u,
u\,f(m_0)+|\lambda|\big)^t\right\|^2&\ge& \left\|\big(-t+u,
u\,f(m_0)\big)^t\right\|^2+|\lambda|^2\\
&\ge&C\,\left(t^2+u^2\right)+|\lambda|^2,
\end{eqnarray*}
where $C$ is a positive constant depending only on $f$. A similar estimate thus holds for
$\nabla_{\theta,\tau}\Phi$, that is,
$$
\big\|\nabla_{\theta,\tau}\Phi\big\|\ge C\,\sqrt{t^2+u^2}+|\lambda|,
$$
for $\lambda\rightarrow -\infty$, possibly for a different constant $C>0$.

The differential operator
$L=:\big\|\nabla_{\theta,\tau}\Phi\big\|^{-2}\cdot
\big(\partial_\theta\Phi\cdot\partial _\theta+\partial_\tau\Phi\cdot\partial_\tau\big)$
satisfies $L\left(e^{i\Phi}\right)=1$; furthermore, its coefficients are rational functions in
$(t,u,\lambda)$ of homogenous degree $-1$, and having homogenous denumerator of
degree $2$ and positive definite in the range $-\lambda,u,t>0$.

Using $L$ to integrate by parts in $d\theta\,d\tau$ in (\ref{eqn:S-microlocal-general}), we conclude
\begin{lem}
Uniformly in $x$, as $\lambda\rightarrow -\infty$
we have
$S_{\chi\,e^{-i\lambda (\cdot)}}(x,x)=O\left(|\lambda|^{-\infty}\right)$.
\end{lem}

Next let us focus on the asymptotics for $\lambda\rightarrow +\infty$, always assuming
$x=x_0+\mathbf{n}$ as above. In this range, we can make the change of variables $t\mapsto \lambda\,t$,
$u\mapsto \lambda\,u$ in  (\ref{eqn:S-microlocal-general}), so as to obtain
\begin{eqnarray}
\label{eqn:S-microlocal-general+infty}
\lefteqn{S_{\chi\,e^{-i\lambda (\cdot)}}(x,x)}\\
&\sim&
\lambda^2\,e^{-i\lambda\,\tau_0}\int_0^{+\infty}\int_0^{+\infty}\left(\int_{-\infty}^{+\infty}\int_{-\pi}^\pi\int_{\mathbb{C}^\mathrm{d}}
e^{i\lambda\Phi_1}\,\mathcal{A}_1\,d\tau\,d\theta\,d\mathbf{v}\right)\,dt\,du,\nonumber
\end{eqnarray}
where now
\begin{equation}
\label{eqn:defn-Phi-A-1}
\begin{array}{ccl}
  \Phi_1 & =: & t\,\psi\big(x,x_0+(\theta,\mathbf{v})\big)
+u\,\psi\left(\phi^X_{-\tau-\tau_0}\big(x_0+(\theta,\mathbf{v})\big),x\right)-\tau,\,\,\,\,
  \mathrm{and}  \\
  \mathcal{A}_1 & =: & a_{\tau+\tau_0}\big(\lambda\,t,x,x_0+(\theta,\mathbf{v})\big)\,
s_{\tau+\tau_0}\big(\lambda\,u,x_0+(\theta,\mathbf{v}),x\big)\\
&&\cdot\chi(\tau+\tau_0)\,\varrho'(\theta,\mathbf{v})\,
\mathcal{V}(\theta,\mathbf{v})
\end{array}
\end{equation}
(and the integrand is supported where $|\tau|<\epsilon$).

Now I claim that only a negligible contribution is lost as $\lambda\rightarrow +\infty$,
if in (\ref{eqn:S-microlocal-general+infty}) integration in $dt\,du$ is restricted to a suitable
compact domain.
Indeed, for some $C\gg 0$ let $\rho\in \mathcal{C}^\infty_0\big((1/2C,2C)\big)$ be $\ge 0$ and $\equiv 1$
on $(1/C,C)$. Then
$$
S_{\chi\,e^{-i\lambda (\cdot)}}(x,x)\sim
S_{\chi\,e^{-i\lambda (\cdot)}}(x,x)_1+
S_{\chi\,e^{-i\lambda (\cdot)}}(x,x)_2,
$$ where $S_{\chi\,e^{-i\lambda (\cdot)}}(x,x)_j$
is defined by the right hand side of (\ref{eqn:S-microlocal-general+infty}), with $A_1$ replaced by
$\mathcal{A}_1'=:\mathcal{A}_1\cdot \rho\big(\|(t,u)\|\big)$ for $j=1$, and by
$\mathcal{A}_1''=:\mathcal{A}_1\cdot \big[1- \rho\big(\|(t,u)\|\big)\big]$ for $j=2$.

\begin{lem}
\label{lem:rapid-decay-u-t}
$S_{\chi\,e^{-i\lambda (\cdot)}}(x,x)_2=O\left(\lambda^{-\infty}\right)$ as
$\lambda\rightarrow +\infty$.
\end{lem}

\begin{proof}
By construction,
$\big(x_0+\mathbf{n},x_0+(\theta,\mathbf{v})\big)$ is at distance $O(\delta)$ from
$(x_0,x_0)$, and $\left(\phi^X_{-\tau-\tau_0}\big(x_0+(\theta,\mathbf{v})\big),x_0+\mathbf{n}\right)$
is at distance $O(\delta+\epsilon)$ from
$(x_0,x_0)$. Since $d_{(y,y)}\psi=(\alpha_y,-\alpha_y)$
for any $y\in X$, using (\ref{eqn:action-in-tau})  we conclude
\begin{equation}
\label{eqn:estimate-on-gradient-theta-tau}
\begin{array}{ccl}
\partial_\theta\Phi_1&=&-t+
u+t\,O\big(\delta\big)+u\,O\big(\delta+\epsilon\big),\\
\partial_\tau\Phi_1&=&u\,f(m_0)-1+u\,O\big(\delta+\epsilon\big).
\end{array}
\end{equation}
On the support of $1-\rho\big(\|(t,u)\|\big)$, we have either $\|(t,u)\|\ge C$,
or else $\|(t,u)\|\le 1/C$.

Where $\|(t,u)\|\ge C$, given (\ref{eqn:estimate-on-gradient-theta-tau}) we have
$\|\nabla_{\theta,\tau}\Phi\|^2\ge C'\,\left(1+t^2+u^2\right)$. Introducing the
operator $L_1=\|\nabla_{\theta,\tau}\Phi\|^{-2}\,\big(\partial_\theta\Phi_1\cdot \partial_\theta
+\partial_\tau\Phi_1\cdot \partial_\tau\big)$, we have
$\lambda^{-1}\cdot L_1\left(e^{i\lambda\Phi_1}\right)=e^{i\lambda\Phi_1}$; arguing as before,
we may use $L_1$ to integrate by parts in $d\theta\,d\tau$ and obtain that the contribution to
the asymptotics of the locus where $\|(t,u)\|\ge C$ is $O\left(\lambda^{-\infty}\right)$ for
$\lambda\rightarrow +\infty$.

On the other hand, where $\|(t,u)\|\le 1/C$ we have $\|\nabla_{\theta,\tau}\Phi\|^2\ge 1/2$;
the same argument implies that the corresponding contribution is also $O\left(\lambda^{-\infty}\right)$.

\end{proof}

By (\ref{eqn:S-microlocal-general+infty}) and Lemma \ref{lem:rapid-decay-u-t}, for some $D=2C\gg 0$ we have
for $\lambda\rightarrow +\infty$
\begin{eqnarray}
\label{eqn:S-microlocal-general+infty-compact}
\lefteqn{S_{\chi\,e^{-i\lambda (\cdot)}}(x,x)}\\
&\sim&
\lambda^2\,e^{-i\lambda\,\tau_0}\int_{\mathbb{C}^\mathrm{d}}\left(\int_{-\pi}^\pi\int_{1/D}^{D}\int_{1/D}^{D}\int_{-\epsilon}^{+\epsilon}
e^{i\lambda\Phi_1}\,\mathcal{A}_1'\,d\theta\,dt\,du\,d\tau\right)\,d\mathbf{v}.\nonumber
\end{eqnarray}

To proceed, we need some asymptotic control on $\mathrm{dist}_X\left(x,\phi^X_{-\tau-\tau_0}(x)\right)$
for $x\rightarrow \mathrm{Fix}\left(\phi^X_{-\tau_0}\right)$ and $|\tau|<\epsilon$.
Let $Z_0$ be the connected component of $\mathrm{Fix}\left(\phi^X_{-\tau_0}\right)$ through $x_0$,
and set $F_0=:\pi(Z_0)$. By assumption, $F_0$ is a compact symplectic submanifold of $M$, and
$Z_0=\pi^{-1}(F_0)$. Let furthermore $M'$ be a tubular neighborhood of $F_0$ and set $X'=:\pi^{-1}(M')$.

Let us define, for $C,\lambda>0$,
\begin{equation}
\label{eqn:V-lambda}
V_\lambda=:\left\{x\in X':
\mathrm{dist}_X(x,Z_0)\ge C\,\lambda^{-7/18}\right\}.
\end{equation}
For $\lambda\rightarrow +\infty$, $U_\lambda=:X\setminus V_\lambda$ is a shrinking tubular neighborhood
of $Z_0$, invariant under the circle action;
in particular it is disjoint from the other connected components of
$\mathrm{Fix}\left(\phi^X_{-\tau_0}\right)$.

Since $\pi$ is a Riemannian submersion, $\mathrm{dist}_X(x,Z_0)=\mathrm{dist}_M(m,F_0)$ if $m=\pi(x)$.
Therefore, if $x\in V_\lambda$ and $m=\pi(y)$ then
$$
\mathrm{dist}_X\left(x,\phi^X_{-\tau_0}(x)\right)\ge \mathrm{dist}_M\left(m,\phi^M_{-\tau_0}(m)\right)
\ge aC\,\lambda^{-7/18}
$$
for some $a>0$.
More generally,
\begin{lem}
\label{lem:distance-tau-n}
Suppose that $\tau_0$ is a very clean period of $\phi^X$.
Then, perhaps after suitably decreasing $\epsilon$, we may find a constant $C'>0$
such that if $\lambda\gg 0$,
$x\in V_\lambda$ and $\tau\in (\tau_0-\epsilon,\tau_0+\epsilon)$ then
$$
\mathrm{dist}_X\left(x,\phi^X_{-\tau}(x)\right)\ge \mathrm{dist}_M\left(m,\phi^M_{-\tau}(m)\right)
\ge C'\,\lambda^{-7/18},
$$
where $m=\pi(x)$.
\end{lem}

This is essentially Lemma 3.2 of \cite{p-ltf}, with some adjustments
to the present more general setting.

\begin{proof}
It suffices to prove the second inequality. Suppose $2\mathrm{r}=\dim F_0$.

Locally in $F_0$, we may find
unitary trivializations of the restricted tangent bundle $\left.TM\right|_{F_0}$, inducing for any
$m'\in F_0$ unitary identifications
$$
T_{m'}M\cong \mathbb{R}^{2\mathrm{r}}\oplus \mathbb{R}^{2(\mathrm{d}-\mathrm{r})},
$$ under which
$\mathbb{R}^{2\mathrm{r}}\oplus (\mathbf{0})$ corresponds to $T_{m'}F_0$, and $(\mathbf{0})\oplus \mathbb{R}^{2(\mathrm{d}-\mathrm{r})}$
to the Riemannian orthocomplement $\mathcal{N}_{m'}=\left (T_{m'}F_0\right)^\perp\subseteq T_{m'}M$
(the latter is generally different from the symplectic normal space $N_{m'}$).

Using the exponential maps for $F_0$ and $M$ as in \cite{p-ijgmmp}, we may then construct a
 smoothly varying family $\xi_{m'}$
 of preferred local coordinates centered at $m'\in M_0$, with the following properties.

For $m'\in F_0$ and $\mathbf{v}\in \mathbb{R}^{2\mathrm{d}}$ sufficiently small, let us set
$$
m'+\mathbf{v}=:\xi_{m'}(\mathbf{v}),
$$
and
let us write the general $\mathbf{v}\in \mathbb{R}^{2\mathrm{d}}$ as $\mathbf{v}=\mathbf{t}\oplus \mathbf{p}$, where
$\mathbf{t}\in \mathbb{R}^{2\mathrm{r}}$ and $\mathbf{p}\in \mathbb{R}^{2(\mathrm{d}-\mathrm{r})}$.
Then
$m'+\mathbf{t}\oplus \mathbf{0}\in F_0$ for every $\mathbf{t}\in \mathbb{R}^{2\mathrm{r}}$; furthermore,
for any $\mathbf{p}\in \mathbb{R}^{2(\mathrm{d}-\mathrm{r})}$
the path
$m'+\mathbf{0}\oplus (s\mathbf{p})$, defined for $s\in (-\epsilon,\epsilon)$,
is perpendicular to $F_0$ at $m'$ for $s=0$.

If $m'\in F_0$, then $\phi^M_{-(\tau+\tau_0)}(m')=\phi^M_{-\tau}(m')\in F_0$ for any $\tau\in \mathbb{R}$;
therefore, by construction for $\tau\sim 0$ we have
$\phi^M_{-(\tau+\tau_0)}(m')=m'+\upsilon(\tau)\oplus \mathbf{0}$, where $\upsilon(\tau)=\tau\,\upsilon_f(m')+O\left(\tau^2\right)$
(here $\upsilon_f(m')\in \mathbb{R}^{2\mathrm{r}}$ in the local coordinates $\xi_{m'}$).

On the other hand, for
$\mathbf{p}\in \mathbb{R}^{2(\mathrm{d}-\mathrm{r})}$ we have
$$
\phi^M_{-\tau_0}\big(m'+\mathbf{0}\oplus \mathbf{p}\big)=
m'+d_{m'}\phi^M_{-\tau_0}(\mathbf{0}\oplus \mathbf{p})+O\left(\|\mathbf{p}\|^2\right),
$$
where again $d_{m'}\phi^M_{-\tau_0}:\mathbb{R}^{2\mathrm{d}}\rightarrow \mathbb{R}^{2\mathrm{d}}$ is computed
in local coordinates.

Since $F_0$ is invariant under $\phi^M_{-\tau_0}$, by construction
$d_{m'}\phi^M_{-\tau_0}\big(\mathbb{R}^{2\mathrm{r}}\oplus (\mathbf{0})\big)=\mathbb{R}^{2\mathrm{r}}\oplus (\mathbf{0})$;
in fact, the latter is the $+1$-eigenspace of $d_{m'}\phi^M_{-\tau_0}$. Therefore,
$d_{m'}\phi^M_{-\tau_0}\left((\mathbf{0})\oplus \mathbb{R}^{2(\mathrm{d}-\mathrm{r})}\right)\subseteq \mathbb{R}^{2\mathrm{d}}$
is transverse to $\mathbb{R}^{2\mathrm{r}}\oplus (\mathbf{0})$. Thus
$$
d_{m'}\phi^M_{-\tau_0}\big(\mathbf{0}\oplus \mathbf{p}\big)=\mathbf{p}'\oplus \mathbf{p}'',
$$
where
$\mathbf{p}'\in \mathbb{R}^{2\mathrm{r}}$ and $\mathbf{p}''\in \mathbb{R}^{2(\mathrm{d}-\mathrm{r})}$
depend linearly on $\mathbf{p}\in \mathbb{R}^{2(\mathrm{d}-\mathrm{r})}$, and $\mathbf{p}''\neq \mathbf{0}$
whenever $\mathbf{p}\neq \mathbf{0}$.

We also have $\mathbf{p}''\neq \mathbf{p}$, for any
$\mathbf{p}\in \mathbb{R}^{2(\mathrm{d}-\mathrm{r})}\setminus \{\mathbf{0}\}$. Indeed, suppose
$\mathbf{p}\in \mathbb{R}^{2(\mathrm{d}-\mathrm{r})}$ and
$d_{m'}\phi^M_{-\tau_0}\big(\mathbf{0}\oplus \mathbf{n}\big)=\mathbf{p}'\oplus \mathbf{p}$
for some $\mathbf{p}'\in \mathbb{R}^{2\mathrm{r}}$. Let
$\mathrm{id}_{m'}:T_{m'}M\rightarrow T_{m'}M$ be the identity.
Then
$$
\left(d_{m'}\phi^M_{-\tau_0}-\mathrm{id}_{m'}\right)\big(\mathbf{0}\oplus \mathbf{p}\big)
=\mathbf{p}'\oplus \mathbf{0}\in \ker\left(d_{m'}\phi^M_{-\tau_0}-\mathrm{id}_{m'}\right).
$$
Since $\tau_0$ is a very clean period, $\mathrm{im}\left(d_{m'}\phi^M_{-\tau_0}-\mathrm{id}_{m'}\right)\cap
\ker\left(d_{m'}\phi^M_{-\tau_0}-\mathrm{id}_{m'}\right)=(\mathbf{0})$. Thus $\mathbf{p}'=\mathbf{0}$,
which implies $\mathbf{0}\oplus \mathbf{p}\in \ker\left(d_{m'}\phi^M_{-\tau_0}-\mathrm{id}_{m'}\right)=T_{m'}F_0$,
whence $\mathbf{p}=\mathbf{0}$.

Summing up, for $\mathbf{p}\sim \mathbf{0}\in \mathbb{R}^{2(\mathrm{d}-\mathrm{r})}$ and $\tau\sim 0$ we have
\begin{eqnarray}
\label{eqn:action-in-local}
\phi^M_{-(\tau+\tau_0)}\big(m'+\mathbf{0}\oplus \mathbf{n}\big)=m'+\big(H(\tau,\mathbf{p}),K(\tau,\mathbf{p})\big),
\end{eqnarray}
where
\begin{eqnarray*}
%\label{eqn:H-K}
 H(\tau,\mathbf{p}) & = & \tau\,\upsilon_f(m)+\mathbf{p}'+O\left(\tau^2+|\tau|\,\|\mathbf{p}\|+
\|\mathbf{p}\|^2\right), \\
 K(\tau,\mathbf{n}) & = & \mathbf{p}''+O\left(|\tau|\,\|\mathbf{p}\|+
\|\mathbf{p}\|^2\right),
\end{eqnarray*}
and there exists $C_1>0$ such that $\left\|\mathbf{p}''-\mathbf{p}\right\|\ge C_1\,\|\mathbf{p}\|$,
for any $m'\in F_0$ and $\mathbf{p}\in \mathbb{R}^{2(\mathrm{d}-\mathrm{r})}$.

Since preferred local coordinates are isometric at the origin, we have
\begin{equation}
\label{eqn:distance-bound-general}
\mathrm{dist}_M(m'+\mathbf{u},m'+\mathbf{w})\ge \frac 12\,\|\mathbf{u}-\mathbf{w}\|,
\end{equation}
for any $\mathbf{u},\mathbf{w}\in \mathbb{R}^{2\mathrm{d}}$ of sufficiently small norm.
Given this, (\ref{eqn:action-in-local}) implies that for
$(\tau,\mathbf{p})\sim (0,\mathbf{0})\in \mathbb{R}\times  \mathbb{R}^{2\mathrm{d}}$
we have
\begin{eqnarray}
\label{eqn:distance-bound}
\mathrm{dist}_M\Big(\phi^M_{-(\tau+\tau_0)}\big(m'+\mathbf{0}\oplus \mathbf{p}\big),
m'+\mathbf{0}\oplus \mathbf{p}\Big)\ge\frac 13\,C_1\,\|\mathbf{p}\|,
\end{eqnarray}
uniformly in $m'\in F_0$.
On the other hand, in the same range we have
$$
\frac 12\,\|\mathbf{p}\|\le \mathrm{dist}_M(m'+\mathbf{0}\oplus \mathbf{n},F_0\big)\le
2\,\|\mathbf{p}\|.
$$
Therefore, any $m\in \pi(V_\lambda)$ may be written in a locally
unique manner as $m=m'+\mathbf{0}\oplus \mathbf{p}$, with $m'\in F_0$ and
$\mathbf{p}\in \mathbb{R}^{2(\mathrm{d}-\mathrm{r})}$ such
that $\|\mathbf{p}\|\ge (C/2)\,\lambda^{-7/18}$. We then see from (\ref{eqn:distance-bound}) that
\begin{eqnarray*}
\mathrm{dist}_M\Big(\phi^M_{-(\tau+\tau_0)}\big(m\big),
m\Big)\ge\frac 16\,C_1\,C\,\lambda^{-7/18}
\end{eqnarray*}
if $|\tau|<\epsilon$.
\end{proof}

We can now prove statement 2 in the Theorem. We shall denote by $C_1,C_2,\ldots$ suitable positive
constants that may be chosen uniformly in the pertinent local data.

By the above, any $x\in V_\lambda$ may be written in a locally unique manner
as $x=x_0+\mathbf{0}\oplus\mathbf{p}$, where $x_0\in Z_0$, $\mathbf{p}\in \mathbb{R}^{2(\mathrm{d}-\mathrm{r})}$,
$\|\mathbf{p}\|\ge C_1\,\lambda^{-7/18}$, and the expression is interpreted in a smoothly varying
Heisenberg local coordinate system centered at $x_0$.

By Lemma \ref{lem:distance-tau-n}, if $m_0=\pi(x_0)$ and $m=\pi(x)=m_0+\mathbf{0}\oplus\mathbf{p}$
then
\begin{eqnarray*}
\mathrm{dist}_X\left(x,\phi^X_{-\tau}(x)\right)
\ge
\mathrm{dist}_M\left(m,\phi^M_{-\tau}(m)\right)
\ge C_2\,\lambda^{-7/18},
\end{eqnarray*}
for any $\tau\in (\tau_0-\epsilon,\tau_0+\epsilon)$.
Given any $\mathbf{v}\in \mathbb{R}^{2\mathrm{d}}$ with $\|\mathbf{v}\|\le \delta$,
therefore,
\begin{eqnarray*}
C_2\,\lambda^{-7/18}\le \mathrm{dist}_M\left(m,m_0+\mathbf{v}\right)+\mathrm{dist}_M\left(m_0+\mathbf{v},
\phi^M_{-\tau}(m)\right),
\end{eqnarray*}
so that
$$
\max\left\{\mathrm{dist}_M\left(m,m_0+\mathbf{v}\right),
\mathrm{dist}_M\left(m_0+\mathbf{v},
\phi^M_{-\tau}(m)\right)\right\}\ge \frac{1}{2}\,C_2\,\lambda^{-7/18}.
$$

Where $\mathrm{dist}_M\left(m,m_0+\mathbf{v}\right)\ge (C_2/2)\,\lambda^{-7/18}$,
for any $\theta\in (-\pi,\pi)$ we also have
$\mathrm{dist}_X\big(x,x_0+(\theta,\mathbf{v})\big)\ge (C_2/2)\,\lambda^{-7/18}$.
Either by the discussion of $\psi$ in \cite{bs} or by a direct inspection using the
computations in \S 3 of \cite{sz}, this implies
\begin{eqnarray*}
\Big|\psi\Big(x,x_0+(\theta,\mathbf{v})\Big)\Big|\ge
\Im\psi\Big(x,x_0+(\theta,\mathbf{v})\Big)\ge C_3\,\,\lambda^{-7/9}.
\end{eqnarray*}
Given this and (\ref{eqn:defn-Phi-A-1}), iterated integration by parts in $t$ in
(\ref{eqn:S-microlocal-general+infty-compact}) introduces at
each step a factor $\lambda^{-2/9}$, and therefore the corresponding contribution to the
asymptotics is $O\left(\lambda^{-\infty}\right)$.

On the other hand, where $\mathrm{dist}_M\left(m_0+\mathbf{v},
\phi^M_{-\tau}(m)\right)\ge (C_2/2)\,\lambda^{-7/18}$
we have $\mathrm{dist}_X\big(x_0+(\theta,\mathbf{v}),
\phi^X_{-\tau}(x)\big)\ge (C_2/2)\,\lambda^{-7/18}$; therefore
\begin{eqnarray*}
\Big|\psi\Big(x_0+(\theta,\mathbf{v}),
\phi^X_{-\tau}(x)\Big)\Big|\ge
\Im\psi\Big(x_0+(\theta,\mathbf{v}),
\phi^X_{-\tau}(x)\Big)\ge C_4\,\,\lambda^{-7/9}.
\end{eqnarray*}
Iteratively integrating by parts in $u$,
we then deduce that also in this case the contribution to the asymptotics is $O\left(\lambda^{-\infty}\right)$.
This proves statement 2 in Theorem \ref{thm:main}.

Let us now consider statement 3. By our geometric assumption, working with smoothly varying Heisenberg
local coordinates centered at points $x'\in Z_0$, any $x\in X$ with $\mathrm{dist}_X(x,Z_0)\le C\,\lambda^{-7/18}$
may be written in a locally unique manner as $x=x'+\mathbf{n}$, where $x'\in Z_0$ and
$\mathbf{n}\in N_{m'}\subseteq T_{m'}M\cong \mathbb{R}^{2\mathrm{r}}\oplus \mathbb{R}^{2(\mathrm{d}-\mathrm{r})}$
satisfies $\|\mathbf{n}\|\le C'\,\lambda^{-7/18}$. In fact, in terms of the previous decomposition
$\mathbf{n}=\mathbf{t}(\mathbf{n})\oplus \mathbf{p}(\mathbf{n})$ where $\mathbf{t}(\mathbf{n})\in \mathbb{R}^{2\mathrm{r}}$ and $\mathbf{p}(\mathbf{n})\in \mathbb{R}^{2(\mathrm{d}-\mathrm{r})}$ depend linearly on $\mathbf{n}$, and
in intrinsic notation
the projection map $\mathbf{n}\mapsto\mathbf{p}(\mathbf{n})$ is a linear isomorphism $N_{m'}\cong \mathcal{N}_{m'}$;
therefore, there exist $a,A>0$ such that $a\,\big\|\mathbf{p}(\mathbf{n})\big\|\le
\|\mathbf{n}\|\le A\,\big\|\mathbf{p}(\mathbf{n})\big\|$.
So to determine the asymptotics on the shrinking tubular neighborhood $U_\lambda=X\setminus V_\lambda$
(see (\ref{eqn:V-lambda}) we may set $x=x_0+\mathbf{n}$, where $x_0\in Z_0$ and $\|\mathbf{n}\|\le C\,\lambda^{-7/18}$
(perhaps after changing the value of the positive constant $C>0$).

We can modify (\ref{eqn:S-microlocal-general+infty-compact}) as follows:

\begin{lem}
There exists $\beta\in \mathcal{C}^\infty_0\left(\mathbb{R}\times \mathbb{R}\times \mathbb{C}^\mathrm{d}\right)$
with $\beta\ge 0$ and $\equiv 1$ in a neighborhood of the origin, such that
the following holds. Set
$$
\beta_\lambda(\tau,\theta,\mathbf{v})=:\beta\left(\lambda^{7/18}\,(\tau,\theta,\mathbf{v})\right).
$$
Then uniformly in $x=x_0+\mathbf{n}\in U_\lambda$ we have
\begin{eqnarray}
\label{eqn:S-microlocal-general+infty-compact-restricted}
\lefteqn{S_{\chi\,e^{-i\lambda (\cdot)}}(x,x)}\\
&\sim&
\lambda^2\,e^{-i\lambda\,\tau_0}\int_{\mathbb{C}^\mathrm{d}}
\left(\int_{-\pi}^\pi\int_{1/D}^{D}\int_{1/D}^{D}\int_{-\epsilon}^{+\epsilon}
e^{i\lambda\Phi_1}\,\mathcal{A}_1''\,d\theta\,dt\,du\,d\tau\right)\,d\mathbf{v},\nonumber
\end{eqnarray}
as $\lambda\rightarrow +\infty$, where
$\mathcal{A}_1''=:\mathcal{A}_1'\cdot \beta_\lambda$.
\end{lem}

\begin{proof}
We are assuming $x=x_0+\mathbf{n}$, with $\mathbf{n}\in N_{m_0}$
(where $m_0=\pi(x_0)$) and $\|\mathbf{n}\|\le C\,\lambda^{-7/18}$.
On the other hand,
since Heisenberg local coordinates are isometric at the origin, we have
\begin{equation}
\label{eqn:distance-bound-general-X}
\mathrm{dist}_X\big(x_0+(\vartheta,\mathbf{u}),x_0+\mathbf{w})\ge
\frac 12\,\big\|(\vartheta,\mathbf{u}-\mathbf{w})\big\|,
\end{equation}
for any sufficiently small
$\vartheta\in \mathbb{R}$ and $\mathbf{u},\mathbf{w}\in \mathbb{R}^{2\mathrm{d}}$.
Therefore, if $\big\|(\theta,\mathbf{v})\big\|\ge 2C\,\lambda^{-7/18}$, say, we have for some constant
$C'>0$
$$
\mathrm{dist}_X\Big(x_0+\mathbf{n},x_0+(\theta,\mathbf{v})\Big)\ge C'\,\lambda^{-7/18},
$$
and so by (\ref{eqn:defn-Phi-A-1})
$$
\big |\partial _t\Phi_1\big |=\left |\psi\Big(x,x_0+(\theta,\mathbf{v})\Big)\right|
\ge C''\,\lambda^{-7/9}
$$
for a suitable constant $C''>0$.
As before, iteratively integrating by parts in $dt$ in (\ref{eqn:S-microlocal-general+infty-compact})
we conclude
that only a rapidly decreasing contribution is lost if integration in
$d\theta\,d\mathbf{v}$ is restricted to an open neighborhood of the origin of radius $2C\,\lambda^{-7/18}$.
In other words perhaps after discarding a negligible contribution we may multiply the amplitude $\mathcal{A}_1'$
in (\ref{eqn:S-microlocal-general+infty-compact}) by a factor of the form
$\beta'_\lambda(\theta,\mathbf{v})=:\beta'\left(\lambda^{7/18}\,\big\|(\theta,\mathbf{v})\big\|\right)$,
where $\beta':\mathbb{R}\rightarrow [0,1]$ is a compactly supported bump function, identically equal to $1$
on $(-2C,2C)$.

In this range,
if in addition $|\tau|\ge (D/2)\,\lambda^{-7/18}$ for some sufficiently large $D>0$ then
by (\ref{eqn:action-in-tau})
$$
\mathrm{dist}_X\Big(\phi^X_{-(\tau+\tau_0)}\big(x_0+(\theta,\mathbf{v})\big),x_0+\mathbf{n}\Big)\ge
C'''\lambda^{-7/18},
$$
since $f(m_0)>0$.
It follows from (\ref{eqn:defn-Phi-A-1}) that
$$
\big |\partial _u\Phi_1\big |=\left |\psi\Big(\phi^X_{-(\tau+\tau_0)}\big(x_0+(\theta,\mathbf{v})\big),x\Big)\right|
\ge D''\,\lambda^{-7/9}.
$$
Integrating by parts in $du$, again we conclude that the corresponding contribution to the asymptotics is
$O\left(\lambda ^{-\infty}\right)$ as $\lambda\rightarrow +\infty$.
Thus, we may further multiply$A_1'\cdot \beta'_\lambda$ by an appropriate bump function in $\tau$
of the form $\beta''_\lambda(\tau)=:\beta''\left(\lambda^{7/18}\tau\right)$, where
$\beta''(s)=\beta'(s/D)$ for some fixed $D\gg 0$. with no effect on the asymptotics.

Finally, we may let $\beta=:\beta'\cdot \beta''$.
\end{proof}

Let us now set $x=x_0+\mathbf{n}/\sqrt{\lambda}$ for a fixed $x_0\in Z_0$ and $\mathbf{n}\in N_{m_0}$, with the understanding that $\|\mathbf{n}\|\le C\,\lambda^{1/9}$, and make the change of integration variables
$(\tau,\theta,\mathbf{v})\mapsto (\tau,\theta,\mathbf{v})/\sqrt{\lambda}$, so that integration in $d\tau\,d\theta\,d\mathbf{v}$
is now over a ball centered at the origin of radius $O\left(\lambda^{1/9}\right)$.
In the following computations, we shall loosely denote by $R_k$ a generic smooth function vanishing to order $k$
at the origin on some appropriate Euclidean space (which may change from line to line).

By \S 3 of \cite{bsz}, with $\psi_2$ as in (\ref{eqn:psi-2}),
\begin{eqnarray}
\label{eqn:psi-1st-term-expanded-1}
\lefteqn{t\,\psi\left(x_0+\frac{\mathbf{n}}{\sqrt{\lambda}},
x_0+\left(\frac{\theta}{\sqrt{\lambda}},\frac{\mathbf{v}}{\sqrt{\lambda}}\right)\right)}\\
&=&it\,\left[1-e^{-i\theta/\sqrt{\lambda}}\right]
-\frac{it}{\lambda}\,\psi_2\big(\mathbf{n},\mathbf{v}\big)\,e^{-i\theta/\sqrt{\lambda}}+
t\,R_3\left(\frac{\mathbf{n}}{\sqrt{\lambda}},\frac{\mathbf{v}}{\sqrt{\lambda}}\right)\,e^{-i\theta/\sqrt{\lambda}},  \nonumber\\
&=&-\frac{t\theta}{\sqrt{\lambda}}+\frac{it}{\lambda}\left[\frac 12\,\theta^2-\psi_2(\mathbf{n},\mathbf{v})\right]+
t\,R_3\left(\frac{\theta}{\sqrt{\lambda}},\frac{\mathbf{n}}{\sqrt{\lambda}},
\frac{\mathbf{v}}{\sqrt{\lambda}}\right),  \nonumber
\end{eqnarray}
where $R_3=O\left(\lambda^{-7/6}\right)$ in the present range.

Similarly, in view of (\ref{eqn:action-in-tau}) we have
\begin{eqnarray}
\label{eqn:psi-2nd-term-expanded-1}
\lefteqn{u\,\psi\Big(\phi^X_{-(\tau+\tau_0)}\big(x_0+(\theta,\mathbf{v})\big),x_0+\mathbf{n}\Big)}\\
&=& iu\,\left[1-e^{iE_\lambda}\right]-\frac{iu}{\lambda}\,\psi_2\big(B_\lambda,\mathbf{n}\big)\,e
^{iE_\lambda}+u\,R_3\left(\frac{B_\lambda}{\sqrt{\lambda}},\frac{\mathbf{n}}{\sqrt{\lambda}}\right)\,
e^{iE_\lambda},  \nonumber
\end{eqnarray}
where
$$
E_\lambda=\frac{1}{\sqrt{\lambda}}\,\big(\theta+\tau\,f(m_0)\big)+
\frac 1\lambda\,\tau\,\omega_{m_0}\big(\upsilon_f(m),A\mathbf{v}\big)+
R_3\left(\frac{\tau}{\sqrt{\lambda}},
\frac{\mathbf{v}}{\sqrt{\lambda}}\right),
$$
and
$$
B_\lambda=\Big(A\mathbf{v}-\tau\,\upsilon_f(m_0)\Big)+
\sqrt{\lambda}\,R_2\left(\frac{\tau}{\sqrt{\lambda}},\frac{\mathbf{v}}{\sqrt{\lambda}}\right).
$$

Expanding the exponential, we get
\begin{eqnarray}
\label{eqn:psi-2nd-term-expanded-2}
\lefteqn{iu\,\left[1-e^{iE_\lambda}\right]
=\frac{u}{\sqrt{\lambda}}\,\big(\theta+\tau\,f(m_0)\big)}\\
&&+
\frac{u}{\lambda}\,\left[\frac i2\,\big(\theta+\tau\,f(m_0)\big)^2+
\tau\,\omega_{m_0}\big(\upsilon_f(m_0),A\mathbf{v}\big)\right]+u\,R_3\left(\frac{\theta}{\sqrt{\lambda}},\frac{\tau}{\sqrt{\lambda}},
\frac{\mathbf{v}}{\sqrt{\lambda}}\right);
%O\left(\lambda^{-7/6}\right),
\nonumber
\end{eqnarray}
furthermore,
\begin{eqnarray}
\label{eqn:psi-2nd-term-expanded-3}
\psi_2\big(B_\lambda,\mathbf{n}\big)=
\psi_2\big(A\mathbf{v}-\tau\,\upsilon_f(m_0),\mathbf{n}\big)+
\lambda\,R_3\left(\frac{\mathbf{v}}{\sqrt{\lambda}},\frac{\tau}{\sqrt{\lambda}},
\frac{\mathbf{n}}{\sqrt{\lambda}}\right).\nonumber
\end{eqnarray}

The analogue of (\ref{eqn:psi-1st-term-expanded-1}) then is
\begin{eqnarray}
\label{eqn:psi-1st-term-expanded-4}
\lefteqn{u\cdot\psi\left(\phi^X_{-(\tau+\tau_0)}
\left(x_0+\left(\frac{\theta}{\sqrt{\lambda}},\frac{\mathbf{v}}{\sqrt{\lambda}}\right)\right),
x_0+\frac{\mathbf{n}}{\sqrt{\lambda}}\right)}\\
&=&
\frac{u}{\sqrt{\lambda}}\Big(\theta+\tau\,f(m_0)\Big)\nonumber\\
&&+\frac{u}{\lambda}\,\left[\frac i2\,\big(\theta+\tau\,f(m_0)\big)^2+
\tau\,\omega_{m_0}\big(\upsilon_f(m_0),A\mathbf{v}\big)
-i\,\psi_2\big(A\mathbf{v}-\tau\,\upsilon_f(m_0),\mathbf{n}\big)\right]\nonumber\\
&&+u\,R_3\left(\frac{\mathbf{n}}{\sqrt{\lambda}},\frac{\tau}{\sqrt{\lambda}},
\frac{\theta}{\sqrt{\lambda}},\frac{\mathbf{v}}{\sqrt{\lambda}}\right),\nonumber
\end{eqnarray}
where the remainder is $O\left(\lambda^{-7/6}\right)$ on the domain of integration.

With the rescaling by $\lambda^{-1/2}$,
(\ref{eqn:S-microlocal-general+infty-compact-restricted}) may be rewritten:
\begin{eqnarray}
\label{eqn:S-microlocal-general+infty-compact-rescaled}
\lefteqn{S_{\chi\,e^{-i\lambda (\cdot)}}
\left(x_0+\frac{\mathbf{n}}{\sqrt{\lambda}},x_0+\frac{\mathbf{n}}{\sqrt{\lambda}}\right)}\\
&\sim&
\lambda^{1-\mathrm{d}}\,e^{-i\lambda\,\tau_0}\int_{\mathbb{C}^\mathrm{d}}
\left(\int_{-\infty}^{+\infty}\int_{1/D}^{D}\int_{1/D}^{D}\int_{-\infty}^\infty
e^{i\lambda\Phi_2}\,\mathcal{A}_2\,d\theta\,dt\,du\,d\tau\right)\,d\mathbf{v},\nonumber
\end{eqnarray}
where $\Phi_2=\Phi_1$ and $\mathcal{A}_2=\mathcal{A}_1''$
with $\mathbf{n},\mathbf{v},\theta,\tau$ replaced by
the rescaled arguments $\mathbf{n}/\sqrt{\lambda},\mathbf{v}/\sqrt{\lambda},\theta/\sqrt{\lambda},\tau/\sqrt{\lambda}$.
Integration in $(\theta,\tau,\mathbf{v})$
is still compactly supported, but on an expanding domain of radius
$O\left(\lambda^{1/9}\right)$.

Given (\ref{eqn:defn-Phi-A-1}), (\ref{eqn:psi-1st-term-expanded-1}) and
(\ref{eqn:psi-1st-term-expanded-4}) we have
\begin{eqnarray}
\label{eqn:fase-espansa}
i\lambda \Phi_2&=&i\sqrt{\lambda}\,\Upsilon + \Theta\\
&&+\lambda\,\left[t\,R_3\left(\frac{\theta}{\sqrt{\lambda}},\frac{\mathbf{n}}{\sqrt{\lambda}},
\frac{\mathbf{v}}{\sqrt{\lambda}}\right)+u\,R_3\left(\frac{\mathbf{n}}{\sqrt{\lambda}},\frac{\tau}{\sqrt{\lambda}},
\frac{\theta}{\sqrt{\lambda}},\frac{\mathbf{v}}{\sqrt{\lambda}}\right)\right],  \nonumber
\end{eqnarray}
where
\begin{equation}
\label{eqn:defn-Upsilon}
  \Upsilon =:\theta (u-t)+\tau\,\big(u\,f(m_0)-1\big), \,\,\,\,\,\mathrm{and}
\end{equation}
\begin{eqnarray}
\label{eqn:defn-Theta}
  \Theta & =: & -\frac t2\,\theta^2-\frac u2\,\big(\theta+\tau\,f(m_0)\big)^2+t\,\psi_2(\mathbf{n},\mathbf{v}) \\
   & & +u\,\psi_2\big(A\mathbf{v}-\tau\,\upsilon_f(m_0),\mathbf{n}\big)+
   i\,u\tau\,\omega_{m_0}\big(\upsilon_f(m_0),A\mathbf{v}\big).\nonumber
\end{eqnarray}

We can then rewrite (\ref{eqn:S-microlocal-general+infty-compact-rescaled}) as follows
\begin{eqnarray}
\label{eqn:S-microlocal-general+infty-compact-rescaled-2}
\lefteqn{S_{\chi\,e^{-i\lambda (\cdot)}}\left(x_0+\frac{\mathbf{n}}{\sqrt{\lambda}},x_0+\frac{\mathbf{n}}{\sqrt{\lambda}}\right)}\\
&\sim&
\lambda^{1-\mathrm{d}}\,e^{-i\lambda\,\tau_0}\int_{\mathbb{C}^\mathrm{d}}
\left(\int_{-\infty}^{+\infty}\int_{1/D}^{D}\int_{1/D}^{D}\int_{-\infty}^\infty
e^{i\sqrt{\lambda}\Upsilon}\,e^{i\Theta}\cdot \mathcal{B}\,d\theta\,dt\,du\,d\tau\right)\,d\mathbf{v},\nonumber
\end{eqnarray}
where
\begin{equation}\label{eqn:defn-B}
\mathcal{B}=:e^{\lambda\,\left[t\,R_3\left(\frac{\theta}{\sqrt{\lambda}},\frac{\mathbf{n}}{\sqrt{\lambda}},
\frac{\mathbf{v}}{\sqrt{\lambda}}\right)+u\,R_3\left(\frac{\mathbf{n}}{\sqrt{\lambda}},\frac{\tau}{\sqrt{\lambda}},
\frac{\theta}{\sqrt{\lambda}},\frac{\mathbf{v}}{\sqrt{\lambda}}\right)\right]}\,\mathcal{A}_2.
\end{equation}

The real part of $\Theta$ is
\begin{eqnarray}
\label{eqn:new-amplitude-real}
\Re(\Theta)&=&-\frac u2\,\big(\theta+\tau\,f(m_0)\big)^2-\frac t2\,\theta^2\\
&&-\frac t2\,\big\|\mathbf{n}-\mathbf{v}\big\|^2-\frac u2\,\big\|A\mathbf{v}-\tau\,\upsilon_f(m_0)-\mathbf{n}\big\|^2.
\nonumber
\end{eqnarray}

\begin{lem}
\label{lem:exponential-bound}
There exists $a>0$ such that
$$
\Re(\Theta)<-a\,\left(\theta^2+\tau^2+\|\mathbf{n}\|^2+\|\mathbf{v}\|^2\right).
$$
for $u,t\in (1/D, D)$.
\end{lem}

\begin{proof}
Clearly,
 \begin{eqnarray}
\label{eqn:new-amplitude-real-1}
\Re(\Theta)&\le& -\frac 1{2D}\,\left[\big(\theta+\tau\,f(m_0)\big)^2+\theta^2\right.\nonumber\\
&&\left.
\,\,\,\,\,\,\,\,\,\,\,\,\,+\big\|\mathbf{n}-\mathbf{v}\big\|^2+\big\|A\mathbf{v}-\tau\,\upsilon_f(m_0)-\mathbf{n}\big\|^2\right].
\end{eqnarray}
With $N=\mathrm{im}\left(A-I\right)\subseteq \mathbb{R}^{2\mathrm{d}}$,
let us consider the linear map
$$F=F_{m_0}:
\mathbb{R}^2\times N\times\mathbb{R}^{2\mathrm{d}}\longrightarrow
\mathbb{R}^2\times \mathbb{R}^{2\mathrm{d}}\times\mathbb{R}^{2\mathrm{d}}
$$
given by
$$
F\big((\theta,\tau,\mathbf{n},\mathbf{v})\big)=\big(\theta,\theta+\tau\,f(m_0),\mathbf{n}-\mathbf{v},
A\mathbf{v}-\tau\,\upsilon_f(m_0)-\mathbf{n}\big).
$$
On $\ker(F)$, since $f(m_0)>0$, we need to have
$\theta=\tau=0$ and $\mathbf{n}=\mathbf{v}=A\mathbf{v}$, whence $\mathbf{n}-A\mathbf{n}=\mathbf{0}$.
But then $\mathbf{n}\in \mathrm{im}(A-I)\cap
\ker(A-I)$, and the latter intersection is the null space by assumption. Thus $\mathbf{n}=\mathbf{v}=\mathbf{0}$
and $f$ is injective. The claim follows from this and (\ref{eqn:new-amplitude-real-1}).

\end{proof}

We may use Taylor expansion of $\mathcal{B}$ in $\big(\theta/\sqrt{\lambda},\tau/\sqrt{\lambda},
\mathbf{v}/\sqrt{\lambda},\mathbf{n}/\sqrt{\lambda}\big)$ at the origin to express it as an asymptotic
expansion in descending powers of $\lambda^{1/2}$, with an $N$-th step remainder
bounded by $\lambda^{2\mathrm{d}-N/2}\,Q_N(\theta,\tau,\mathbf{v},\mathbf{n})$, for appropriate $r,N>0$,
and a certain polynomial $Q_N$.
% (with coefficients depending on $u$ and $t$).
By Lemma \ref{lem:exponential-bound}, this expansion may be integrated term
by term to obtain an asymptotic expansion for (\ref{eqn:S-microlocal-general+infty-compact-rescaled-2}).
More explicitly, we can write
$
\mathcal{B}=\mathcal{B}^{(N)}+\mathcal{R}^{(N)},
$
where
$\left|\mathcal{R}^{(N)}\right|\le C\,\lambda^{2\mathrm{d}-(N+1)/2}\,Q_N(\theta,\tau,\mathbf{v},\mathbf{n})$,
while, recalling that $\mathcal{V}=1/(2\pi)$ at the origin, we have
\begin{equation}
\label{eqn:dominant-term}
\mathcal{B}^{(N)}
=\frac{\rho_{\tau_0}(x_0)}{2\pi}\,\left(\frac \lambda\pi\right)^{2\mathrm{d}}\,t^\mathrm{d}\,u^{\mathrm{d}}
\left(1+\sum_{j=1}^N\lambda^{-j/2}P_j(u,t,\theta,\tau,\mathbf{n},\mathbf{v})\right),
\end{equation}
for appropriate polynomials $P_j$. We deduce from this and (\ref{eqn:S-microlocal-general+infty-compact-rescaled-2})
that
\begin{eqnarray}
\label{eqn:smoothed-composition-rescaled-Nth-1}
S_{\chi\,e^{-i\lambda (\cdot)}}\left(x_0+\frac{\mathbf{n}}{\sqrt{\lambda}},
x_0+\frac{\mathbf{n}}{\sqrt{\lambda}}\right)
\sim
\int_{\mathbb{C}^{\mathrm{d}}}I_\lambda(\mathbf{n},\mathbf{v})\,d\mathbf{v}
+O\left(\lambda^{1+\mathrm{d}-(N+1)/2}\right),\nonumber
\end{eqnarray}
where
\begin{eqnarray}
\label{eqn:smoothed-composition-rescaled-Nth-2}
I_\lambda(\mathbf{n},\mathbf{v})&=:&
\frac{\rho_{\tau_0}(x_0)}{2\pi^\mathrm{d}}\left(\frac\lambda\pi\right)^{1+\mathrm{d}}\,e^{-i\lambda\tau_0}\\
&&\cdot
\int_{-\infty}^{+\infty}\int_{1/D}^{D}\int_{1/D}^{D}
\int_{-\infty}^{\infty}
e^{i\sqrt{\lambda}\Upsilon}\,e^{\Theta}\cdot S_\lambda^{(N)}\,d\tau\,dt\,du\,d\theta,  \nonumber
%&&+O\left(\lambda^{1+\mathrm{d}-(N+1)/2}\right),\nonumber
\end{eqnarray}
with
\begin{equation}
\label{eqn:s-n-lambda}
S_\lambda^{(N)}=:t^\mathrm{d}\,u^{\mathrm{d}}
\left(1+\sum_{j=1}^N\lambda^{-j/2}P_j(u,t,\theta,\tau,\mathbf{n},\mathbf{v})\right).
\end{equation}

Let us evaluate asymptotically the inner integral (\ref{eqn:smoothed-composition-rescaled-Nth-1}),
by viewing it as an oscillatory integral in $\sqrt{\lambda}$ with phase $\Upsilon$.
To begin with, from (\ref{eqn:defn-Upsilon}) we obtain
$\partial_t\Upsilon=-\theta$, $\partial_u\Upsilon=\theta+\tau\,f(m_0)$. Since
$f(m_0)>0$, this implies
$\left\|\nabla_{t,u}\Upsilon\right\|\ge C\sqrt{\theta^2+\tau^2}$ for some $C>0$.
We may use the partial differential operator
$R=:\left\|\nabla_{t,u}\Upsilon\right\|^{-2}\,
\Big(\partial _t\Upsilon\cdot \partial_t+\partial _u\Upsilon\cdot \partial_u\Big)$ to integrate by parts
in $t,u$, and obtain that for any $a>0$
the contribution to the asymptotics of the locus where $|(\theta,\tau)|\ge a$
is rapidly decreasing; since on the other hand the domain of integration in $d\mathbf{v}$
is a ball of radius $O\left(\lambda^{1/9}\right)$, the same holds for the contribution to
the global integral (equivalently, one may invoke again
the rapidly decreasing exponential in
$\mathbf{v}$). Therefore the asymptotics are unchanged, if the integrand is multiplied by
a compactly supported bump function in $(\theta,\tau)$, identically $=1$ near the origin.
We shall leave this bump function implicit in the following, and proceed assuming that
the integrand is compactly supported.

We have:

\begin{lem}
\label{lem:stationary-point-upsilon}
The phase $\Upsilon=\Upsilon(\theta,t,u,\tau)$ has the unique stationary point
$$
(\theta_0,t_0,u_0,\tau_0)=\left(0,\frac{1}{f(m_0)},\frac{1}{f(m_0)},0\right).
$$
The stationary point is non-degenerate, and the Hessian matrix there is
$$
H(\Upsilon)_0=\left(
              \begin{array}{cccc}
                0 & -1 & 1 & 0 \\
                -1 & 0 & 0 & 0 \\
                1 & 0 & 0 & f(m_0) \\
                0 & 0 & f(m_0) & 0 \\
              \end{array}
            \right),
$$
with determinant $\det \big(H(\Upsilon)_0\big)=f(m_0)^2$. In particular,
the quadratic form associated to $H(\Upsilon)_0$ has vanishing signature.
\end{lem}

We can then apply the stationary phase Lemma in $\sqrt{\lambda}$. We have
\begin{equation}
\label{eqn:hessian-lambda}
\sqrt{\det\left(\frac{\sqrt{\lambda}}{2\pi i}\cdot H(\Upsilon)_0\right)}
=\frac{\lambda}{(2\pi)^2}\cdot f(m_0).
\end{equation}

Applying Theorem 7.7.5 of \cite{h-libro}, we obtain
\begin{eqnarray}
\label{eqn:leading-term}
I_\lambda(\mathbf{n},\mathbf{v})&\sim&2\pi\cdot \frac{\rho_{\tau_0}(x_0)}{\pi^{\mathrm{d}}}\cdot
\,\left(\frac \lambda\pi\right)^{\mathrm{d}}\frac{e^{-i\lambda\tau_0}}{f(m_0)^{2\mathrm{d}+1}}\\
&&\cdot e^{f(m_0)^{-1}\big[\psi_2(\mathbf{n},\mathbf{v})+\psi_2(A\mathbf{v},\mathbf{n})\big]}\cdot \left(1+\sum_{j\ge 1}\lambda^{-j/2}P_j(\mathbf{n},\mathbf{v})\right),\nonumber
\end{eqnarray}
with the $P_j$'s polynomials.
The $N$-th step remainder is uniformly bounded by
$\lambda^{\mathrm{d}-(N+1)/2}\,e^{-a(\|\mathbf{n}\|^2+\|\mathbf{v}\|^2)}$ for some $a>0$
(and integration in $d\mathbf{v}$
is over a ball of radius $O\left(\lambda^{1/9}\right)$), so the expansion may be integrated term by term
in $d\mathbf{v}$.
We end up with an asymptotic expansion of the form
\begin{eqnarray}
\label{eqn:asymptotic-expansion-j}
S_{\chi\,e^{-i\lambda (\cdot)}}
\left(x_0+\frac{\mathbf{n}}{\sqrt{\lambda}},x_0+\frac{\mathbf{n}}{\sqrt{\lambda}}\right)\sim
\sum_{j\ge 0}\lambda^{\mathrm{d}-j/2}p_j,
\end{eqnarray}
where, with $P_0=1$, the $j$-th coefficient is given by
\begin{eqnarray}
\label{eqn:jth-coefficient}
p_j&=&2\pi\cdot\frac{\rho_{\tau_0}(x_0)}{\pi^{2\mathrm{d}}}\cdot\frac{e^{-i\lambda\tau_0}}{f(m_0)^{2\mathrm{d}+1}}\\
&&\cdot\int_{\mathbb{C}^\mathrm{d}}e^{f(m_0)^{-1}\big[\psi_2(\mathbf{n},\mathbf{v})+\psi_2(A\mathbf{v},\mathbf{n})\big]}
P_j(\mathbf{n},\mathbf{v})\,d\mathbf{v}.\nonumber
\end{eqnarray}

The standard identification $\mathbb{C}^\mathrm{d}\cong \mathbb{R}^{2\mathrm{d}}$ associates
$(\mathbf{x}^t\,\mathbf{y}^t)^t\in \mathbb{R}^{2\mathrm{d}}$ to $\mathbf{x}+i\mathbf{y}\in \mathbb{C}^\mathrm{d}$, for any
$\mathbf{x},\mathbf{y}\in \mathbb{R}^d$. Multiplication by $i$ then corresponds to the $2\mathrm{d}\times 2\mathrm{d}$ matrix
$$
J_0=\left(
      \begin{array}{cc}
        0 & -I_\mathrm{d} \\
        I_\mathrm{d} & 0 \\
      \end{array}
    \right),
$$
and the standard symplectic structure $\omega_0$ is represented by $-J_0$.
Using this and the change of variables $\mathbf{v}=\mathbf{u}+\mathbf{n}$, we obtain
\begin{eqnarray}
\label{eqn:exponent}
\lefteqn{\psi_2(\mathbf{n},\mathbf{v})+\psi_2(A\mathbf{v},\mathbf{n})}\\
%&=&\psi_2\left(\mathbf{n}',\mathbf{v}'\right)+\psi_2\left(A\mathbf{v}',\mathbf{n}'\right)\\
&=&-i\,\omega_0\left(\mathbf{n},\mathbf{v}\right)-\frac 12\,\left\|\mathbf{n}-\mathbf{v}\right\|^2
-i\,\omega_0\left(A\mathbf{v},\mathbf{n}\right)-\frac 12\,\left\|A\mathbf{v}-\mathbf{n}\right\|^2
\nonumber\\
&=&\psi_2\left(A\mathbf{n},\mathbf{n}\right)
+i\,\omega_0\left(A^{-1}\mathbf{n}-\mathbf{n},\mathbf{u}\right)-\frac 12 \left(\|\mathbf{u}\|^2+\|A\mathbf{u}\|^2\right)
-g_0\big(A\mathbf{u},A\mathbf{n}-\mathbf{n}\big)\nonumber\\
&=&\psi_2\left(A\mathbf{n},\mathbf{n}\right)
+i\,\mathbf{u}^tJ_0\left(A^{-1}-I\right)\mathbf{n}-\frac 12 \mathbf{u}^t\big(I+A^tA\big)\mathbf{u}
-\mathbf{u}^tA^t(A-I)\mathbf{n},\nonumber
\end{eqnarray}
where $g_0$ denotes the standard Euclidean structure on $\mathbb{R}^{2\mathrm{d}}$.

As in the Introduction,
let us set $F=:J_0\left(A^{-1}-I\right)=(A^t-I)\,J_0$, $Q=I+A^tA$, $G=:A^t\,(A-I)$; furthermore,
let us set
$\mathbf{n}'=:\mathbf{n}/\sqrt{f(m_0)}$. Then, applying (\ref{eqn:exponent})
and the further change of integration
variable $\mathbf{u}= \sqrt{f(m_0)}\,\mathbf{w}$,
(\ref{eqn:jth-coefficient}) may be rewritten as follows:
\begin{eqnarray}
\label{eqn:jth-coefficient-int}
p_j&=&2\pi\cdot\frac{\rho_{\tau_0}(x_0)}{\pi^{2\mathrm{d}}}\cdot
\frac{1}{f(m_0)^{\mathrm{d}+1}}\cdot e^{-i\lambda\tau_0+
f(m_0)^{-1}\psi_2\left(A\mathbf{n},\mathbf{n}\right)}\\
&&\cdot\int_{\mathbb{C}^\mathrm{d}}
e^{i\,\mathbf{w}^tF\mathbf{n}'-\frac 12 \mathbf{w}^tQ\mathbf{w}
-\mathbf{w}^tG\mathbf{n}'}
Q_j(\mathbf{n}',\mathbf{w})\,d\mathbf{w}\nonumber
\end{eqnarray}
for certain polynomials $Q_j$.

Now let us set
$\mathbf{w}=\mathbf{s}-Q^{-1}G\mathbf{n}'$, and define
$$
\Gamma(\mathbf{n}')=:-i\,{\mathbf{n}'}^{t}\,G^tQ^{-1}F\,\mathbf{n}'+\frac 12\,{\mathbf{n}'}^{t}\,G^tQ^{-1}G\,\mathbf{n}'.
$$
Then for appropriate polynomials $\widetilde{Q}_j$
we can rewrite the integral in
(\ref{eqn:jth-coefficient-int}) as
\begin{eqnarray}
\label{eqn:jth-coefficient-int-1}
e^{\Gamma(\mathbf{n}')}\int_{\mathbb{C}^\mathrm{d}}
e^{i\,\mathbf{s}^tF\mathbf{n}'-\frac 12 \mathbf{s}^tQ\mathbf{s}
}
\widetilde{Q}_j(\mathbf{n}',\mathbf{s})\,d\mathbf{s}.
\end{eqnarray}

Let us consider the leading term of the expansion (\ref{eqn:asymptotic-expansion-j}). Since
\begin{eqnarray}
\label{eqn:basic-integral}
\int_{\mathbb{C}^\mathrm{d}}
e^{i\,\mathbf{s}^tF\mathbf{n}'-\frac 12 \mathbf{s}^tQ\mathbf{s}}\,d\mathbf{s}=
\frac{(2\pi)^{\mathrm{d}}}{\sqrt{\det(Q)}}\,
e^{-\frac 12{\mathbf{n}'}^tF^tQ^{-1}F\mathbf{n}'},
\end{eqnarray}
we obtain from (\ref{eqn:jth-coefficient-int}), (\ref{eqn:jth-coefficient-int-1}) and (\ref{eqn:basic-integral})
that the leading term of the expansion is
\begin{eqnarray}
\label{eqn:jth-coefficient-int-j0}
\lambda^{\mathrm{d}}\,p_0&=&\left(\frac \lambda\pi\right)^\mathrm{d}\cdot\rho_\tau(x_0)\,
\frac{2\pi\,e^{-i\lambda\tau_0}}{f(m_0)^{\mathrm{d}+1}}\cdot\frac{2^\mathrm{d}}{\sqrt{\det(Q)}}
\\
&&\cdot \exp\left(
\frac{1}{f(m_0)}\left[\psi_2(A\mathbf{n},\mathbf{n})
+\Gamma(\mathbf{n})-\frac 12\,\mathbf{n}^tF^tQ^{-1}F\mathbf{n}\right]\right).\nonumber
\end{eqnarray}
By Lemma \ref{lem:RA=PA}, the exponent on the second line of (\ref{eqn:jth-coefficient-int-j0})
is
$$f(m_0)^{-1}\mathbf{n}^t\mathfrak{R}_A\mathbf{n}=
f(m_0)^{-1}\mathbf{n}^t\mathfrak{P}_A\mathbf{n}=\Psi_2^A(\mathbf{n}).$$

Now let us consider the lower order terms of (\ref{eqn:asymptotic-expansion-j}), given by
(\ref{eqn:jth-coefficient-int}) and (\ref{eqn:jth-coefficient-int-1}).
Up to a multiple of $2\pi$, (\ref{eqn:jth-coefficient-int-1}) is $e^{\Gamma(\mathbf{n}')}$
times the evaluation at $-F\mathbf{n}'$ of the Fourier transform of
$e^{-\frac 12{\mathbf{s}}^tQ\mathbf{s}}\cdot \widetilde{Q}_j(\mathbf{n},\mathbf{s})$. The latter Fourier transform,
as a function of $\xi\in \mathbb{R}^{2\mathrm{d}}$, has the
form
$S_j(\mathbf{n},D_\mathbf{\xi})\left(e^{-\frac 12{\xi}^tQ^{-1}\xi}\right)$, where
$S_j(\mathbf{n},D_\xi)$
is an appropriate differential polynomial in $\xi$, whose coefficients are polynomials in $\mathbf{n}$.
It thus still has the form $S_j'(\mathbf{n},\xi)\cdot e^{-\frac 12{\xi}^tQ^{-1}\xi}$
for some polynomial $S_j'$.
Statement 3 of the theorem follows by setting $\xi=-F\mathbf{n}'$.

Finally let us prove the last statement of the Theorem.
Since the asymptotic expansions for the amplitudes $a_\tau$ and $s_\tau$ in (\ref{eqn:R-tau-FIO})
and (\ref{eqn:Pi-tau-FIO}) go down by integer steps, the appearance of half-integer powers of
$\lambda$ in (\ref{eqn:leading-term}) is due to the combined effect of first Taylor expanding the amplitudes
(with various remainders incorporated) in $(\mathbf{v},\mathbf{n},\theta,\tau)/\sqrt{\lambda}$ at the origin
in (\ref{eqn:dominant-term}), and then
applying the stationary phase Lemma in $\sqrt{\lambda}$ in (\ref{eqn:leading-term}).

On the one hand, therefore, $S_\lambda^{(N)}$ in (\ref{eqn:s-n-lambda})
may be decomposed as the sum of terms of the form $\lambda^{-(r_1+r_2)-(d_\mathbf{v}+d_\mathbf{n}+d_\theta+d_\tau)/2}\,
F_{\underline{d}}(\mathbf{v},\mathbf{n},\theta,\tau)$, where $r_j$ is an integer, and $F_{\underline{d}}$ is a polyhomogeneous polynomial of
multidegree $\underline{d}=(d_\mathbf{v},d_\mathbf{n},d_\theta,d_\tau)$.

On the other hand, the inverse of the Hessian matrix $H(\Upsilon)$ in Lemma \ref{lem:stationary-point-upsilon}
is
$$
H(\Upsilon)^{-1}=\frac{1}{f(m_0)}\,\left(
                                     \begin{array}{cccc}
                                       0 & -f(m_0) & 0 & 0 \\
                                       -f(m_0) & 0 & 0 & 1 \\
                                       0 & 0 & 0 & 1 \\
                                       0 & 1 & 1 & 0 \\
                                     \end{array}
                                   \right),
$$
and therefore when apply the stationary phase in $\sqrt{\lambda}$ to each of the previous summands
we introduce (up to scalar multiples irrelevant to the present argument) terms of the form
$$
\left.2\,\lambda^{1-(r_1+r_2)-(d_\mathbf{v}+d_\mathbf{n}+d_\theta+d_\tau)/2-k/2}\,\left(-f(m_0)\cdot\frac{\partial^2}{\partial \theta\partial t}
+\frac{\partial^2}{\partial \tau\partial t}+\frac{\partial^2}{\partial \tau\partial u}\right)^kF_j\right|_{\theta=\tau=0}.
$$
In turn, this splits as a linear combination of terms of the form
$$
\left.\lambda^{1-(r_1+r_2)-(d_\mathbf{v}+d_\mathbf{n}+d_\theta+d_\tau)/2-k/2}\,D^{(k)}_{t,u}\circ \frac{\partial^a}{\partial \theta^a}\circ \frac{\partial^b}{\partial \tau^b}F_j\right|_{\theta=\tau=0,\,t=u=1/f(m_0)},
$$
where $a+b=k$ and $D^{(k)}_{t,u}$ is a differential operator of degree $k$ in $(t,u)$.
Since we are evaluating at the origin in $(\theta,\tau)$, a non-zero contribution is obtained only for $a=d_\theta$, $b=d_\tau$,
so that $k=d_\theta+d_\tau$.
As in (\ref{eqn:jth-coefficient}), we conclude that,
up to factors that we may presently omit, each of the previous terms gives rise in the final expansion to a summand of the form
\begin{equation}
\label{eqn:general-summand-degree}
\lambda^{r-(d_\mathbf{v}+d_\mathbf{n})/2}
\int_{\mathbb{C}^\mathrm{d}}e^{f(m_0)^{-1}\big[\psi_2(\mathbf{n},\mathbf{v})+\psi_2(A\mathbf{v},\mathbf{n})\big]}
S_j(\mathbf{n},\mathbf{v})\,d\mathbf{v},
\end{equation}
where $r$ is an integer, and $S_j$ is bihomogeneous
in $(\mathbf{n},\mathbf{v})$ of bidegree $(d_{\mathbf{n}},d_\mathbf{v})$.

Furthermore, under the change of variables
$\mathbf{v}=\mathbf{u}+\mathbf{n}$ in (\ref{eqn:exponent}),
the rescaling by $\sqrt{f(m_0)}$ in (\ref{eqn:jth-coefficient-int}), and the further change of variables
$\mathbf{w}=\mathbf{s}-Q^{-1}G\mathbf{n}'$ in (\ref{eqn:jth-coefficient-int-1}), the collective degree
$d_\mathbf{v}+d_\mathbf{n}$ is unchanged. Thus, as in (\ref{eqn:jth-coefficient-int-1}),
 (\ref{eqn:general-summand-degree}) splits as
a linear combination of terms of the form
\begin{equation}
\label{eqn:general-summand-degree-1}
\lambda^{r-(d_\mathbf{v}+d_\mathbf{n})/2}
\int_{\mathbb{C}^\mathrm{d}}
e^{i\,\mathbf{s}^tF\mathbf{n}'-\frac 12 \mathbf{s}^tQ\mathbf{s}
}
R_j(\mathbf{n}',\mathbf{s})\,d\mathbf{s},
\end{equation}
where $R_j$ is bihomogenous in $(\mathbf{n},\mathbf{s})$, of
bidegree $(d_{\mathbf{n}},d_{\mathbf{v}})$.

Taking the Fourier transform, this is the evaluation  at $\xi=-F\mathbf{n}'$ of an expression of the form
\begin{equation}
\label{eqn:general-summand-degree-2}
\lambda^{r-(d_\mathbf{v}+d_\mathbf{n})/2}
e^{-\frac 12 \xi^tQ^{-1}\xi
}
S_j(\mathbf{n}',\xi)\,d\mathbf{s},
\end{equation}
where $S_j$ is the sum of various monomials, each of which
has degree $d_\mathbf{n}$ in $\mathbf{n}$, and $d_\xi=d_\mathbf{v}-2k$ for some $k\ge 0$ in $\xi$.
If we set $d'=:d_\xi+d_\mathbf{n}$,
we obtain from each such summand a contribution of the form
\begin{eqnarray}
\label{eqn:general-summand-degree-3}
\lefteqn{\lambda^{r-k-(d_\xi+d_\mathbf{n})/2}
e^{-\frac 12 {\mathbf{n}'}^tF^tQ^{-1}F\mathbf{n}'
}
\widetilde{S}_j(\mathbf{n}',F\mathbf{n}')}\\
&=&\lambda^{r-k-d'/2}
e^{-\frac 12 {\mathbf{n}'}^tF^tQ^{-1}F\mathbf{n}'
}
\widetilde{S}_j'(\mathbf{n}'),\nonumber
\end{eqnarray}
where $r$ and $k$ are integers, and $\widetilde{S}_j'(\mathbf{n}')=:\widetilde{S}_j(\mathbf{n}',F\mathbf{n}')$
is homogeneous of degree $d'$. Thus the exponent of $\lambda$ is integer (respectively, fractional) if and only
if $d'$ is even (respectively, odd).

The proof of the Theorem is complete.

\end{document}